\documentclass[11pt]{article}
\usepackage{amsmath,amssymb,enumerate,amsfonts,mathrsfs,theorem}
\usepackage[a4paper]{geometry}
\pdfoutput=1 
\usepackage{graphicx}
\usepackage{subfigure}
\usepackage{color}
\usepackage{relsize}
\usepackage[retainorgcmds]{IEEEtrantools}
\usepackage{color}
\usepackage[framemethod=TikZ]{mdframed}
\usepackage{lipsum}
\usepackage{soul}

\usepackage[font={footnotesize}]{caption}

\usepackage{fancybox}
\cornersize{0}

\usepackage{tikz}

\theorembodyfont{\it}
\newtheorem{theorem}{Theorem}

\newtheorem{lemma}{Lemma}

\theorembodyfont{\rm}

\definecolor{cooperative}{rgb}{0.78,0.86,0.94}
\definecolor{competitive}{rgb}{0.68,0.72,0.78}
\usepackage{framed}
\usepackage{lipsum}

\def\dotquote{\begin{quote} }
\def\enddotquote{''\end{quote}}

\title{Modeling the modulation of neuronal bursting:\\a singularity theory approach\thanks{This paper presents research results of the Belgian Network
DYSCO (Dynamical Systems, Control, and Optimization),
funded by the Interuniversity Attraction Poles
Programme, initiated by the Belgian State, Science Policy
Office. The scientific responsibility rests with its authors.}}

\author{
Alessio Franci$^{1,*}$,
Guillaume Drion$^{2,3}$, 
\& Rodolphe Sepulchre$^{2,4}$\\
\small $^1$ INRIA Lille-Nord Europe, Orchestron project, 40 avenue Halley F 59650\\ \small Villeneuve d'Ascq, France\\
\small $^2$ Department of Electrical Engineering and Computer Science, University of Liege, Liege, Belgium\\
\small $^3$ Laboratory of Neurophysiology, GIGA Neurosciences, University of Liege, Liege, Belgium\\
\small $^4$ University of Cambridge, Department of Engineering, Trumpington Street,\\ \small Cambridge CB2 1PZ, United Kingdom\\
\small $^*$ Corresponding author. Email: {\tt alessio.franci@inria.fr}
}

\date{}

\begin{document}

\maketitle

\abstract{Exploiting the specific structure of neuron conductance-based models, the paper investigates the mathematical modeling of neuronal bursting modulation. The proposed approach combines singularity theory and geometric singular perturbations to {capture} the geometry of multiple time-scales attractors in the neighborhood of high-codimension singularities. We detect a three-time scale bursting attractor in the universal unfolding of the winged cusp singularity and discuss the physiological relevance of the bifurcation and unfolding parameters in determining a physiological modulation of bursting. The results suggest generality and simplicity in the organizing role of the winged cusp singularity for the global dynamics of conductance based models.}

\section{Introduction}

Bursting is an important signaling  component of neurons, characterized by a periodic {alternation} of bursts and quiescent periods. Bursts are {transient}, but high-frequency trains of spikes, contrasting with the absence of spikes during the quiescent periods. Bursting activity has been recorded in many neurons, both in vitro and in vivo, and electrophysiological recordings show a great variety of bursting time {series}. All neuronal bursters share nevertheless a sharp separation between three different time scales: a fast time-scale for the spike generation, a slow time-scale for the intraburst spike frequency, and an ultra slow time-scale for the inter burst frequency.   Many neuronal models exhibit bursting in some parameter range and many bursting models have been analyzed through bifurcation theory  but the exact {mechanisms modulating neuronal} bursting are still poorly understood, both mathematically and physiologically. In particular, modeling the route to burst, that is the {physiologically observed modulation} from a regular pacemaking activity to a bursting activity, has remained elusive to date. 
{Also many efforts have been devoted at classifying different types of bursters} \cite{Rinzel1987b,Bertram1995,Golubitsky2001,Izhikevich2000}. {But the mathematical mechanisms that allow a same neuron to be modulated across different types are rarely studied, despite their physiological role in homeostatic cell regulation and development} \cite{Liu1998}.

As an attempt to advance the mathematical understanding of neuronal bursting, the present paper exploits the particular structure of conductance based neuronal models to {address} with a local analysis tool the global structure of {bursting} attractors. Rooted in the seminal work of Hodgkin and Huxley \cite{Hodgkin1952}, conductance-based models are nonlinear RC circuits consisting of one capacitance (modeling the cell membrane) in parallel with possibly many voltage sources with voltage dependent conductance (each modeling a specific ionic current). The variables of the model are the membrane potential ($V$) and the gating (activation and inactivation) variables that model the kinetics of each ion channel. The vast diversity of ion channels involved in a particular neuron type leads to high-dimensional models, but all conductance-based models share two central structural assumptions:
\begin{itemize}
\item[] {\it (i)} a classification of gating variables in three well separated time-scales (fast variables - in the range of the membrane potential time scale $\sim 1ms$; slow variables - $5$ to $10$ times slower; and ultra-slow variables - $10$ to hundreds time slower), which roughly correspond to the three time scales of neuronal bursting.
\item[] {\it (ii)} each voltage regulated gating variable $x$ obeys the first-order monotone dynamics $\tau_x(V)\dot x=-x+x_\infty(V)$, which implies that, at steady state, every voltage regulated gating variable is an explicit monotone function of the membrane potential, that is, $x=x_\infty(V)$.
\end{itemize}

Our analysis of neuronal bursting rests on these two structural assumptions. Assumption {\it (i)} suggests a three-time scale singularly perturbed bursting model, whose singular limit provides the skeleton of the bursting attractor. Assumption {\it (ii)} implies that the equilibria of arbitrary conductance-based models are determined by Kirchoff's law (currents sum to zero in the circuit), which provides a single algebraic equation in the sole scalar variable $V$. This remarkable feature calls for singularity theory \cite{Golubitsky1985} to understand the equilibrium structure of the model.

{The results of jointly exploiting timescale separation and singularity theory for neuronal bursting modeling provide the following specific contributions:

$-$The universal unfolding of the winged-cusp singularity is shown to organize a three time-scale burster. The three level hierarchy of singularity theory dictates the hierarchy of timescales: the state variable of the bifurcation problem is the fast variable, the bifurcation parameter is the slow variable, and unfolding parameter(s) are the ultra-slow variable(s). Because the geometric construction is grounded in the algebraic and timescale structure of conductance-based models, the proposed model can be related to detailed conductance- based models through mathematical reduction. We provide general conditions for this mathematical model to be a normal form reduction of an arbitrary conductance-based model. Both the bifurcation parameter and the  unfolding parameters have a clear physiological interpretation.

$-$ The bifurcation parameter is directly linked to the balance between restorative and regenerative slow ion channels, the importance of which was recently studied by the authors in} \cite{Franci2013}. {The modulation of the bifurcation parameter in the proposed three-time scale model provides a geometrically and physiologically meaningful transition from slow tonic spiking to bursting. This ``route to bursting" is known to play a significant role in central nervous system activity} \cite{Sherman2001,Viemari2006,Beurrier1999}. {Its mathematical modeling appears to be novel.

$-$ The three unfolding parameters  modulate in an even slower scale the fast-slow phase portrait of the three-time scale burster.  The affine parameter plays the classical role of an adaptation current that hysterically modulates the slow-fast phase portrait across a parameter range where a stable resting state and a stable spiking limit cycle coexist, thereby creating the bursting attractor. The two remaining unfolding parameters can modulate the bursting attractor across a continuum of bursting types. As a result, transition between differenting bursting waveforms, observed for instance in developing neurons} \cite{Liu1998}, {are geometrically captured as paths in the unfolding space of the winged cusp. The physiological interpretation of this modulation is a straightforward consequence of the clear physiological interpretation of each unfolding parameter.}

The existence of three-time scale bursters in the abstract unfolding of a winged cusp is presented in Section 2. Section 3 focuses on a minimal reduced model of neuronal bursting and uses the insight of singularity theory to describe a physiological route to bursting in this model. Section 4 shows how to trace the same geometry in arbitrary conductance based models. Section \ref{SEC: why the winged cusp} discuss in a less technical way the relevance of the winged-cusp singularity for the modeling of bursting modulation. The technical details of mathematical proofs are presented in an appendix.

\section{Universal unfolding and multi-time scale attractors}

\subsection{A primer on singularity theory}

We introduce here some notation and terminology that will be used extensively in the paper. The interested reader is referred to the main results of Chapters I-IV in \cite{Golubitsky1985} for a comprehensive exposition of the singularity theory used in this paper.

Singularity theory studies scalar bifurcation problems of the form
\begin{equation}\label{EQ: generic bif problem}
g(x,\lambda)=0,\quad x,\lambda\in\mathbb R,
\end{equation}
where $g$ is a smooth function. The variable $x$ denotes the state and $\lambda$ is the bifurcation parameter. The set of pairs $(x,\lambda)$ satisfying (\ref{EQ: generic bif problem}) is called the {\it bifurcation diagram}. {\it Singular points} satisfy $g(x^\star,\lambda^\star)=\frac{\partial g}{\partial x}(x^\star,\lambda^\star)=0$. Indeed, if $\frac{\partial g}{\partial x}(x^\star,\lambda^\star)\neq0$, then the implicit function theorem applies and the bifurcation diagram is necessarily regular at $(x^\star,\lambda^\star)$.

Except for the fold $x^2\pm\lambda=0$, bifurcations are not generic, that is they do not persist under small perturbations. Singularity theory is a robust bifurcation theory: it aims at classifying all possible persistent bifurcation diagrams that can be obtained by small perturbations of a given singularity.

A {\it universal unfolding} of $g(x,\lambda)$ is a parametrized family of functions $G(x,\lambda; \alpha)$, where $\alpha$ lies in the unfolding parameter space $\mathbb R^k$, such that
\begin{itemize}
\item[] 1) $G(x,\lambda;0)=g(x,\lambda)$
\item[] 2) Given any $p(x)$ and a small $\mu>0$, one can find an $\alpha$ near the origin such that the two bifurcation problems $G(x,\lambda;\alpha)=0$ and $g(x,\lambda)+\mu p(x)=0$ are qualitatively equivalent.
\item[] 3) $k$ is the minimum number of unfolding parameters needed to reproduce all perturbed bifurcation diagrams of $g(x,\lambda)$. $k$ is called the codimension of $g(x,\lambda)$.
\end{itemize}
Unfolding parameters are not bifurcation parameters. Instead, they change the qualitative bifurcation diagram of the perturbed bifurcation problem $G(x,\lambda;\alpha)=0$. That is why $\lambda$ is a distinguished parameter in the theory. Historically, this parameter was associated to a slow time, whose evolution lets the dynamics visit the bifurcation diagram in a quasi-steady state manner. It will play the same role in the present paper, where we only consider two singularities and their universal unfolding:\\
the codimension 1 hysteresis
\begin{equation}\label{EQ: stable hy}
g_{hy}^s(x,\lambda)=-x^3-\lambda,
\end{equation}
whose universal unfolding is shown to be \cite[Chapter IV]{Golubitsky1985}
\begin{equation}\label{EQ: stable hy unfold}
G_{hy}^s(x,\lambda;\ \beta)=-x^3-\lambda+\beta x,
\end{equation}
the codimension 3 winged cusp
\begin{equation}\label{EQ: stable cusp}
g_{wcusp}^s(x,\lambda)=-x^3-\lambda^2,
\end{equation}
whose universal unfolding is shown to be \cite[Section III.8 and Chapter IV]{Golubitsky1985}
\begin{equation}\label{EQ: stable cusp unfold}
G_{wcusp}^s(x,\lambda;\ \alpha,\beta,\gamma)=-x^3-\lambda^2+\beta x -\gamma\lambda x -\alpha.
\end{equation}
%In the remainder of the paper, we use the compact form ``w-cusp" instead of ``winged cusp".

The universal unfolding of codimension$\geq$1 bifurcations contains some codimension 1 bifurcation. For instance, the universal unfolding of the winged cusp possesses hysteresis bifurcations on the unfolding parameter hypersurface defined by $\alpha\gamma^2+\beta=0$, $\alpha\leq0$. Even though such bifurcation diagrams are not persistent, they define {\it transition varieties} that separate equivalence classes of persistent bifurcation diagrams, hence, providing a complete classification of persistent bifurcation diagrams.

An unperturbed bifurcation problem assumes the suggestive role of {\it organizing center}: all the perturbed bifurcation diagrams are determined and organized by the unperturbed bifurcation diagram, {which constitutes the most singular situation}. Via the inspection of local algebraic conditions at the singularity, an organizing center provides a {quasi-global} description of all possible perturbed bifurcation diagrams.

\subsection{The hysteresis singularity and spiking oscillations}
\label{SSEC: rel osci in hyst}

\begin{figure}[h!]
\center
\includegraphics[width=0.9\textwidth]{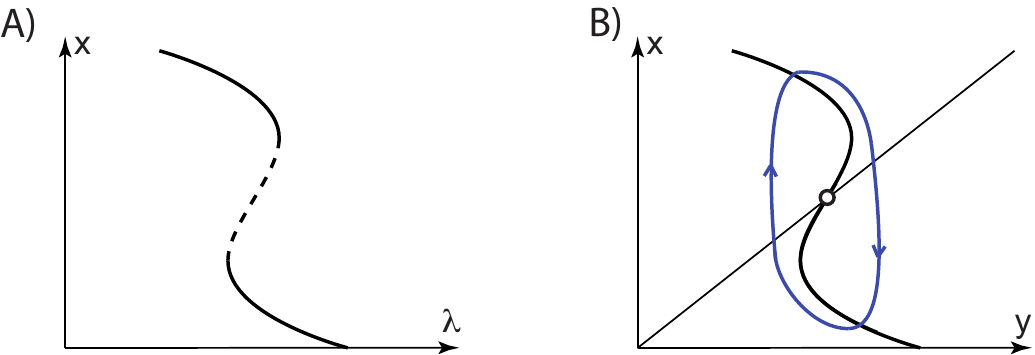}
\caption{{\bf Relaxation oscillations in the universal unfolding of the hysteresis bifurcation.} {\bf A.} A persistent bifurcation diagram of the hysteresis singularity. Branches of stable (resp. unstable) fixed points are depicted as full (resp. dashed)  lines. {\bf B.} Through a slow adaptation of the bifurcation parameter, the bifurcation diagram in A. is transformed into the phase plane of a two-dimensional dynamical system, which still defines a universal unfolding of the hysteresis singularity. The thick full line is the fast subsystem nullcline. The thin full line is the slow subsystem nullcline. The circle denotes an unstable fixed point. For small $\beta>0$, the model exhibits exponentially stable relaxation oscillations (depicted in light blue).}\label{FIG 1}
\end{figure}

The hysteresis singularity has a universal unfolding $-x^3-\lambda+\beta x$ with persistent bifurcation diagram plotted in Figure \ref{FIG 1}A for $\beta>0$. We use this algebraic curve to generate the phase portrait in Fig. \ref{FIG 1}B of the two-time scale model
\begin{IEEEeqnarray}{rCl}\label{EQ: hy dynamics}
\dot x&=&G_{hy}^s(x,\lambda+y;\ \beta)\IEEEyessubnumber\\
&=&-x^3+\beta x-\lambda-y\nonumber\\
\dot y&=&\varepsilon(x-y)\IEEEyessubnumber
\end{IEEEeqnarray}
\noindent Because $y$ is a slow variable, it acts as a slowly varying modulation of the bifurcation parameter in the fast dynamics (\ref{EQ: hy dynamics}a). As a consequence, the global analysis of system (\ref{EQ: hy dynamics}) reduces to a quasi-steady state bifurcation analysis of (\ref{EQ: hy dynamics}a), hence the relationship between Fig. \ref{FIG 1}A and Figure \ref{FIG 1}B.

%{The same would not hold by modulating a non-distinguished (unfolding) parameter, which points to the relevance of singularity theory to our approach. Keeping the bifurcation parameter in} (\ref{EQ: hy dynamics}a) {has a twofold role. First, it ensures that the steady-state structure of} (\ref{EQ: hy dynamics}) {is again organized by a hysteresis singularity, with $\lambda$ as bifurcation parameter and $\beta$ as unfolding parameter. Second, it allows to explore different dynamical regime of the two-dimensional model by changing the base point of the local bifurcation structure. This latter fact will be crucial in the next section, when we will repeat the same construction for the winged-cusp singularity.}

{The following (well known) theorem characterizes a global attractor of} (\ref{EQ: hy dynamics}){, that is the existence of Van-der-Pol type relaxation oscillations in the universal unfolding of the hysteresis.}

\begin{theorem}\label{THM: hy rel oscil}{\bf \cite{Grasman1987},\cite{Mishchenko1980},\cite{Krupa2001c}}
For $\lambda=0$ and for all $0<\beta<1$, there exists $\bar\varepsilon>0$ such that, for all $\varepsilon\in(0,\bar\varepsilon]$, the dynamical system (\ref{EQ: hy dynamics}) possesses an exponentially stable relaxation limit cycle, which attracts all solutions except the equilibrium at $(0,0)$.
\end{theorem}
The familiar reader will recognize in (\ref{EQ: hy dynamics}) a famous model of neurodynamics introduced by FitzHugh \cite{FitzHugh1961}. It is the prototypical planar reduction of spiking oscillations. There is therefore a close relationship between the hysteresis singularity and spike generation.

{It is worth emphasizing that the relationship between singularity theory} (Fig. \ref{FIG 1}A) {and the two-time scale phase portrait} (Fig. \ref{FIG 1}B) {imposes choosing the bifurcation parameter, not an unfolding parameter, as the slow variable. It should also be observed that the slow variable is a deviation from the unfolding parameter $\lambda$ rather than the bifurcation parameter itself. Keeping $\lambda$ as the bifurcation parameter of the two-dimensional dynamics} (\ref{EQ: hy dynamics}) {allows to shape its equilibrium structure accordingly to the universal unfolding of the organizing singularity, in this case, the hysteresis, and will play an important role in the next section.}

\subsection{The winged cusp singularity and rest-spike bistability}
\label{SSEC: rest-spike in wcusp}

We repeat the elementary construction of Section \ref{SSEC: rel osci in hyst} for the codimension-3 winged cusp singularity $-x^3-\lambda^2$. It differs from the hysteresis singularity in the \emph{non-monotonicity} of $g(x,\lambda)$ in the bifurcation parameter, that is $\frac{\partial (-x^3-\lambda^2)}{\partial\lambda}=-2\lambda$ changes sign at the singularity.

\begin{figure}[h!]
\center
\includegraphics[width=0.9\textwidth]{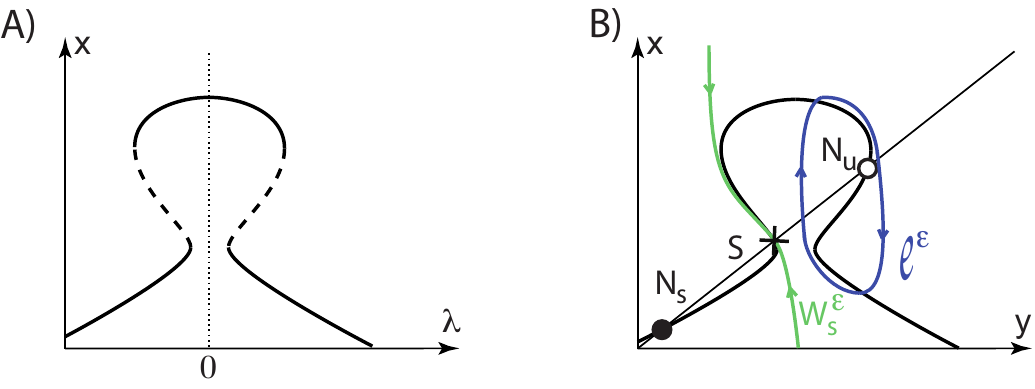}
\caption{{\bf Singularly perturbed rest-spike bistability in the universal unfolding of the winged cusp.} {\bf A.} {\it Mirrored hysteresis} persistent bifurcation diagram of the winged cusp for $\beta>0$, $\alpha<-2\left(\frac{\beta}{3}\right)^{3/2}$, and $\gamma=0$. {\bf B.} A phase plane of (\ref{EQ: cusp dynamics}).}\label{FIG 2}
\end{figure}

Figure \ref{FIG 2}A illustrates an important persistent bifurcation diagram in the unfolding of the winged cusp, obtained for $\gamma=0$, $\beta>0$, and $\alpha<-2\left(\frac{\beta}{3}\right)^{3/2}$. We call it the {\it mirrored hysteresis} bifurcation diagram. The right part ($\lambda>0$) of this bifurcation diagram is essentially the persistent bifurcation diagram of the hysteresis singularity in Figure \ref{FIG 1}A. In that region, $\frac{\partial G_{wcusp}^s}{\partial \lambda}<0$. The left part ($\lambda>0$) is the mirror of the hysteresis and, in that region, $\frac{\partial G_{wcusp}^s}{\partial \lambda}>0$. For $\gamma\neq 0$, the mirroring effect is not perfect, but the qualitative analysis does not change. The hysteresis and its mirror collide in a transcritical singularity for $\alpha=-2\left(\frac{\beta}{3}\right)^{3/2}$. This singularity belongs to the transcritical bifurcation transition variety in the winged cusp unfolding  (see Appendix \ref{SEC: trans variety}). The transcritical bifurcation variety plays an important role in the forthcoming analysis.

We use the algebraic curve in Figure \ref{FIG 2}A to generate the phase portrait in Figure \ref{FIG 2}B of the two-dimensional model
\begin{IEEEeqnarray}{rCl}\label{EQ: cusp dynamics}
\dot x&=&G_{wcusp}^s(x,\lambda+y;\ \alpha,\beta,\gamma)\IEEEyessubnumber\\
&=&-x^3+\beta x-(\lambda+y)^2-\gamma(\lambda+y)x -\alpha\nonumber\\
\dot y&=&\varepsilon(x-y).\IEEEyessubnumber
\end{IEEEeqnarray}
Its fixed point equation
\begin{equation}\label{EQ: cusp dyn fp}
F(x,\lambda,\alpha,\beta,\gamma):=-x^3+\beta x-(\lambda+x)^2-\gamma(\lambda+x)x-\alpha.
\end{equation}
is easily shown to be again a universal unfolding of the winged cusp around $x_{wcusp}:=\frac{1}{3}$, $\lambda_{wcusp}:=0$, $\alpha_{wcusp}:=-\frac{1}{27}$, $\beta_{wcusp}:=-\frac{1}{3}$, $\gamma_{wcusp}:=-2$. The face portrait in Fig. \ref{FIG 2}B is a prototype phase portrait of rest-spike bistability: a stable fixed point coexists with a stable relaxation limit cycle.

\begin{figure}[h!]
\center
\includegraphics[width=0.75\textwidth]{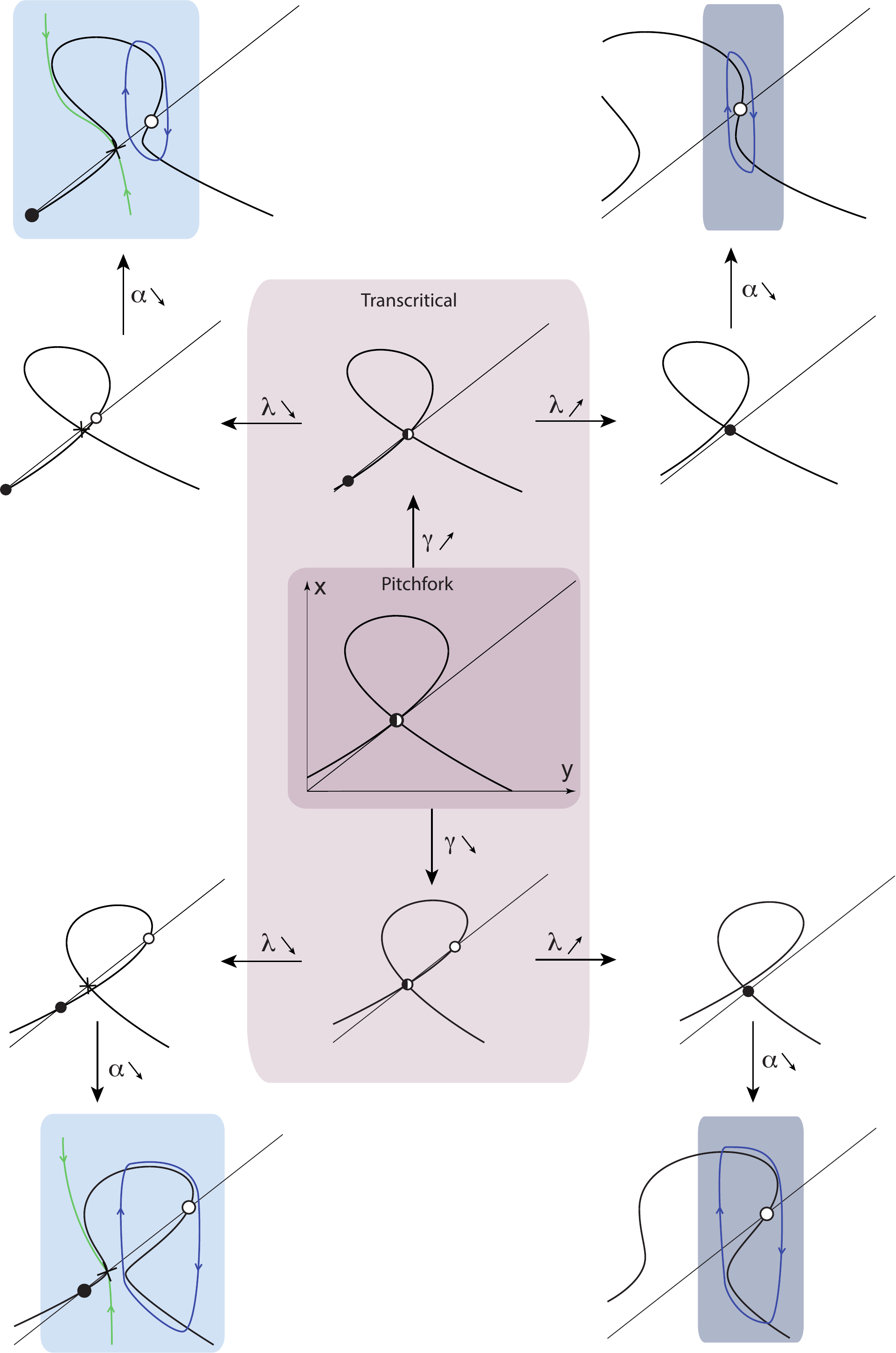}
\caption{{\bf An unfolding of the pitchfork bifurcation variety in (\ref{EQ: cusp dynamics}).}
%The phase portrait of Figure \ref{FIG 1} appears in the top right and bottom right corners. The phase portrait of Fig. \ref{FIG 2} appears in the top left and bottom left corners.
The phase portraits in Figs. \ref{FIG 1} and \ref{FIG 2} both belong to the unfolding of the pitchfork singularity in center. A smooth deformation of the phase portrait of Fig. \ref{FIG 1} into the phase portrait of Fig. \ref{FIG 2} involves a transcritical bifurcation, which degenerate into a pitchfork for a particular value of the unfolding parameter $\gamma$}.\label{FIG 3}
\end{figure}

Similarly to the previous section, the analysis of the singularly perturbed model (\ref{EQ: cusp dynamics}) is completely characterized by the bifurcation diagram of Figure \ref{FIG 2}A. This bifurcation diagram provides a skeleton for the rest-spike bistable phase portrait in Figure \ref{FIG 2}B, as stated in the following theorem. Its proof is provided in Section \ref{SSEC: cusp bist proof}.

\begin{theorem}\label{THM: cusp bist}
For all $\beta>\beta_{wcusp}$, there exist open sets of bifurcation ($\lambda$) and unfolding ($\alpha,\gamma$) parameters near the pitchfork singularity at $(\lambda,\alpha,\gamma)=(\lambda_{PF}(\beta),\alpha_{PF}(\beta),\gamma_{PF}(\beta))$, in which, for sufficiently small $\varepsilon>0$, model (\ref{EQ: cusp dynamics}) exhibits the coexistence of an exponentially stable fixed point $N_s$ and an exponentially stable spiking limit cycle $\ell^\varepsilon$. Their basins of attraction are separated by the stable manifold $W_s^\varepsilon$ of a hyperbolic saddle $S$ {\rm (see Fig. \ref{FIG 2}B)}.
\end{theorem}

Figure \ref{FIG 3} shows the transition in (\ref{EQ: cusp dynamics}) from the hysteresis phase portrait in Figure \ref{FIG 1}B to the bistable phase portrait in Fig. \ref{FIG 2}B through a transcritical bifurcation. Both phase portraits are generated by unfolding the degenerate portrait in Fig. \ref{FIG 3}, center, which belongs to the pitchfork bifurcation variety $(\alpha,\gamma)=(\alpha_{PF}(\beta),\gamma_{PF}(\beta))$, $\beta>\beta_{wcusp}$ (see Appendix \ref{SEC: trans variety}). The transcritical bifurcation variety $\alpha=\alpha_{TC}(\beta,\gamma)$ is obtained  through variations of the unfolding parameter $\gamma$ away from the pitchfork variety. It provides the two phase portraits in Fig. \ref{FIG 3}, center top and bottom. By increasing or decreasing the bifurcation parameter $\lambda$ and decreasing the unfolding parameter $\alpha$ out of the transcritical bifurcation variety, these phase portraits perturb to the generic phase portraits in the corner, corresponding to the qualitative phase portraits in Figures \ref{FIG 1}B and Fig. \ref{FIG 2}B, respectively. The reader of \cite{Franci2012} will recognize the same organizing role of the pitchfork in a planar model of neuronal excitability.

\subsection{A three-time scale bursting attractor in the winged cusp unfolding}
\label{SSEC: cusp burst}

The coexistence of a stable resting state and stable spiking oscillation, or {\it singularly perturbed rest-spike bistability}, makes (\ref{EQ: cusp dynamics}) a good candidate as the slow-fast subsystem of a three-time scale minimal bursting model:
\begin{IEEEeqnarray}{rCl}\label{EQ: 3D cusp dynamics}
\dot x&=&G_{wcusp}^s(x,\lambda+y;\ \alpha+z,\beta,\gamma)\IEEEyessubnumber\\
&=&-x^3+\beta x-(\lambda+y)^2-\gamma(\lambda+y)x -\alpha-z\nonumber\\
\dot y&=&\varepsilon_1(x-y)\IEEEyessubnumber\\
\dot z&=&\varepsilon_2 (-z + a x+b y + c ),\IEEEyessubnumber
\end{IEEEeqnarray}
where {$0<\varepsilon_2\ll\varepsilon_1\ll 1$ and $a,b,c\in\mathbb R$}. The {$z$-dynamics} models the ultra-slow adaptation of the affine unfolding parameter $\alpha$, in such a way that the global attractor of (\ref{EQ: 3D cusp dynamics}) will be determined by {a quasi-static modulation} of (\ref{EQ: 3D cusp dynamics}a) through different persistent bifurcation diagrams.

{Here, again, the role of singularity theory in distinguishing bifurcation and unfolding parameters is crucial. The hierarchy between these parameters and the state variable, formalized in the theory in} \cite[Definition III.1.1]{Golubitsky1985}, {is reflected here in the hierarchy of timescales.}
%$Z$ is a smooth function, and we have chosen to model the ultra-slow adaptation variable through the affine unfolding parameter  $\alpha$. By ultra-slowly modulating the unfolding parameter $\alpha$, the evolution of $z$ smoothly transforms the bifurcation diagram of  (\ref{EQ: 3D cusp dynamics}a).

\begin{figure}[h!]
\center
\includegraphics[width=0.8\textwidth]{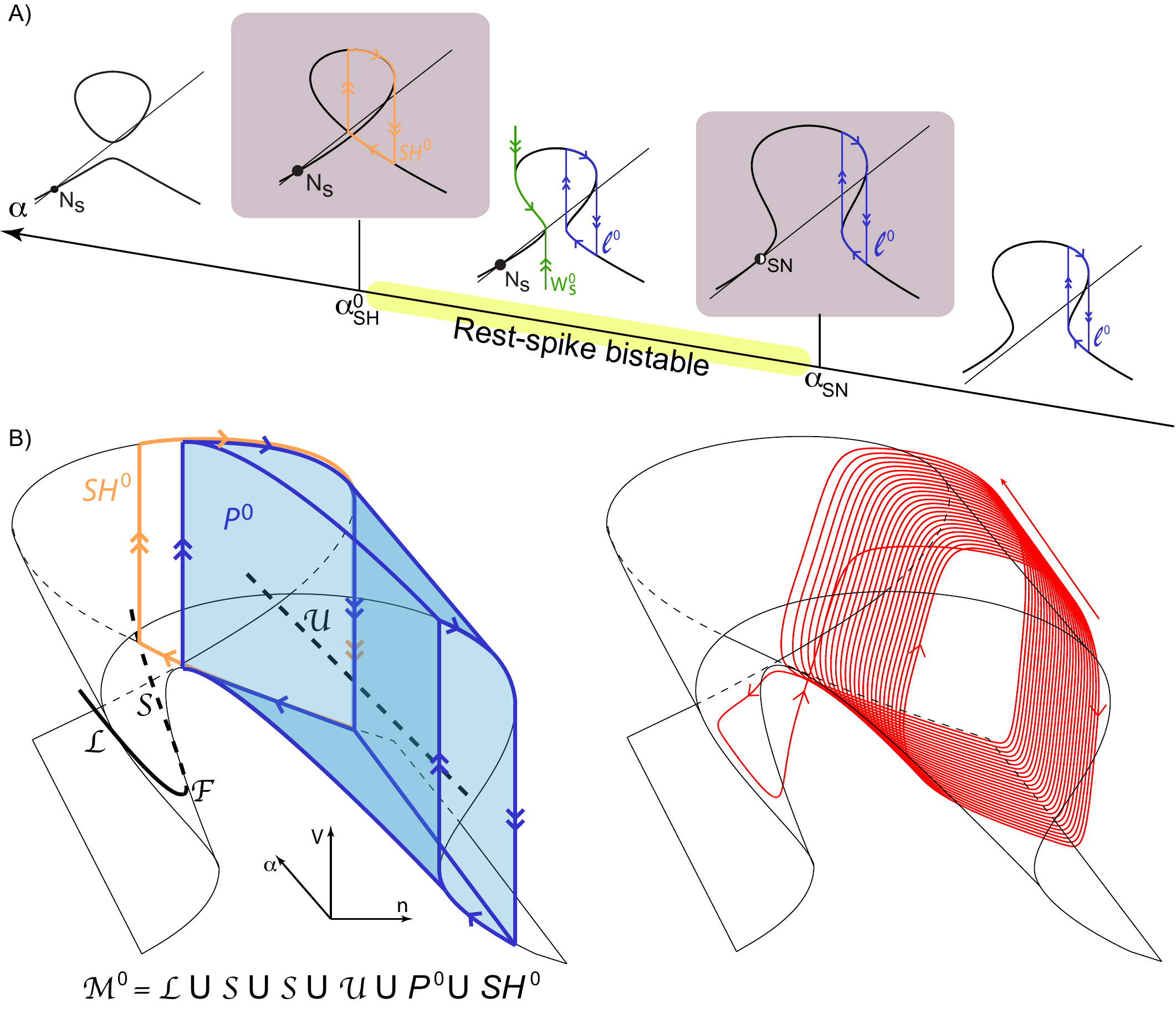}
\caption{{\bf Singular steady-state behavior of (\ref{EQ: cusp dynamics}) through a variation of the unfolding parameter $\alpha$.} {\bf A.} Singular phase portraits of (\ref{EQ: cusp dynamics}) for $\gamma=0$, $\beta=\frac{1}{3}$, and small negative $\lambda$. {\bf B.} Gluing the different invariant sets in {\bf A} leads to the three-dimensional singular invariant set $\mathcal M^0$ (left), which provides a skeleton for a three-time scale bursting attractor (right) in the singularly perturbed system (\ref{EQ: 3D cusp dynamics}). The branch of stable fixed points (resp. saddle point) for $\alpha<\alpha_{SN}$ is drawn as the black solid curve $\mathcal L$ (resp. the black dashed curve $\mathcal S$). The saddle node bifurcation connecting them is denoted by $\mathcal F$. The branch of unstable fixed points is drawn as the black dashed line $\mathcal U$. The branch of stable singular periodic orbit for $\alpha<\alpha_{SH}^0$ is drawn as the blue cilindric surface $P^0$. The singular saddle homoclinic trajectory is drawn as the orange oriented curve $SH^0$.}\label{FIG 4}
\end{figure}

The time scale separation between (\ref{EQ: 3D cusp dynamics}a-\ref{EQ: 3D cusp dynamics}b) and (\ref{EQ: 3D cusp dynamics}-c) makes it possible once again to derive a global analysis of model (\ref{EQ: 3D cusp dynamics}) from the analysis of the steady state behavior of (\ref{EQ: cusp dynamics}) as $\alpha$ is varied. Such  analysis can easily be derived geometrically in the singular limit $\varepsilon_1=0$. It is sketched in Figure \ref{FIG 4}. {For $\alpha\in(\alpha_{SN},\alpha_{SH}^0)$, the singularly perturbed model} (\ref{EQ: cusp dynamics}) {exhibits rest-spike bistability, that is, the coexistence of a stable node $N_s$, a singular stable periodic orbit $\ell^0$, and a singular saddle separatrix $W_s^0$. At $\alpha=\alpha_{SH}^0=-2\left(\frac{\beta}{3}\right)^{3/2}$ the left and right branches of the mirrored hysteresis bifurcation collide in a transcritical singularity that serves as a connecting point for a singular homoclinic trajectory $SH^0$. For $\alpha>\alpha_{SH}^0$, the only (singular) attractor is the stable node $N_s$. At $\alpha=\alpha_{SN}$, the saddle and the stable node merge in a saddle-node bifurcation $SN$. For $\alpha<\alpha_{SN}$, the only attractor is the singular periodic orbit $\ell^0$.} The different singular invariant sets in Figure \ref{FIG 4}A, can be glued together to construct the three-dimensional singular invariant set $\mathcal M_0$ in Figure \ref{FIG 4}B-left.

The singular invariant set $\mathcal M_0$ provides a skeleton for a three-time scale bursting attractor that shadows the branch $\mathcal L$ of stable fixed points in alternation with the branch $P^0$ of (singular) stable periodic orbits, as depicted in Figure \ref{FIG 4}B-right. To prove the existence of such an attractor, we only need to understand how $\mathcal M^0$ perturbs for $\varepsilon_1>0$.

\begin{figure}[h!]
\center
\includegraphics[width=0.5\textwidth]{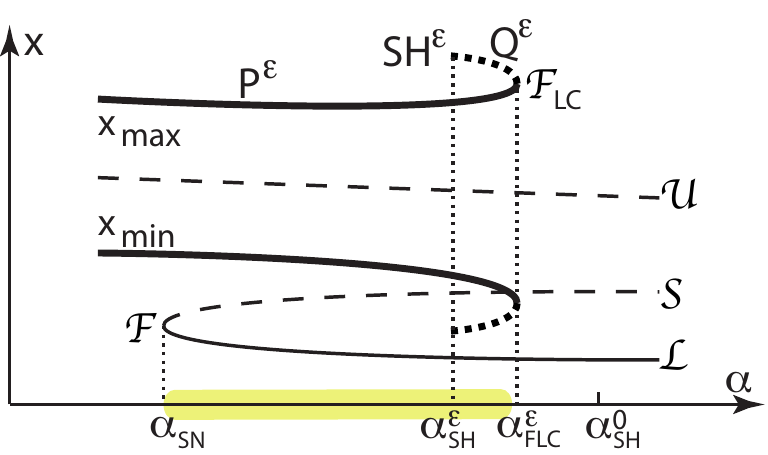}
\caption{{\bf Bifurcation diagram of (\ref{EQ: cusp dynamics}) with respect the unfolding parameter $\alpha$ for sufficiently small $\varepsilon$.} The branch of stable fixed point is depicted as the full thin line $\mathcal{L}$, the branch of saddle point as the dashed thin line $\mathcal S$, and the branch of unstable fixed point as the dashed thin line $\mathcal U$. The branch of stable periodic orbit is depicted by the thick full lines $P^\varepsilon$ and the branch of unstable periodic orbits by the thick dashed lines $Q^\varepsilon$. $SH^\varepsilon$: saddle-homoclinic bifurcation. $\mathcal F_{LC}$: fold limit cicle bifurcation. $\mathcal F$: fold (saddle-node) bifurcation. The yellow strip between the saddle-node and fold limit cicle bifurcation denotes the rest-spike bistable range.}\label{FIG 5}
\end{figure}

Near the singular limit, the branch of singular periodic orbits $P^0$ perturbs to a nearby branch of exponentially stable periodic orbits $P^\varepsilon$ (see Fig. \ref{FIG 5}), whereas the singular homoclinic trajectory $SH^0$ perturbs to an unstable homoclinic trajectory $SH^\epsilon$ (at $\alpha=\alpha_{SH}^\epsilon$). The branch of unstable periodic orbits $Q^\varepsilon$ generated at $SH^\varepsilon$ eventually merges with $P^\varepsilon$ at a fold limit cycle bifurcation $\mathcal F_{LC}$ for some $\alpha_{FLC}^\varepsilon\in(\alpha_{SH}^\epsilon,\alpha_{SH}^0)$. In the whole range $(\alpha_{SN},\alpha_{FLC}^\varepsilon)$, model (\ref{EQ: cusp dynamics}) exhibits the coexistence of a stable fixed point and a stable spiking limit cycle. The details of this analysis are contained in Lemma \ref{LEM: cusp SH} in Section \ref{SSUB: burst att proof}.

We follow \cite{Terman1991,Su2004} to derive conditions on the bifurcation and unfolding parameters in (\ref{EQ: 3D cusp dynamics}a-\ref{EQ: 3D cusp dynamics}b) and
{to place the hyperplane $\dot z=0$ (through a suitable choice of the parameters $a,b,c\in\mathbb R$) such that an ultra-slow variation of} $z$  can hysteretically modulate the slow-fast subsystem (\ref{EQ: 3D cusp dynamics}a-\ref{EQ: 3D cusp dynamics}b) across its bistable range $(\alpha_{SN},\alpha_{FLC}^\varepsilon)$ to obtain stable bursting oscillations. The existence of such bursting oscillations is stated in the following theorem. Its proof is provided in Section \ref{SSUB: burst att proof}.

\begin{theorem}\label{THM: cusp burst}
For all $\beta>\beta_{wcusp}$, there exists an open set of bifurcation ($\lambda$) and unfolding ($\alpha,\gamma$) parameters near the pitchfork singularity at $(\lambda,\alpha,\gamma)=(\lambda_{PF}(\beta),\alpha_{PF}(\beta),\gamma_{PF}(\beta))$ {such that, for all $\lambda,\alpha,\gamma$ in those sets, there exist $a,b,c,\in\mathbb R$ such that,} for sufficiently small $\varepsilon_1\gg\varepsilon_2>0$, model (\ref{EQ: 3D cusp dynamics}) has a hyperbolic bursting attractor.
\end{theorem}

{Theorem} \ref{THM: cusp burst} {uses the two regenerative phase portraits in Fig.} \ref{FIG 3} {left to construct a} {\it bursting} {attractor by modulating the unfolding parameter $\alpha$. The bursting attractor directly rests upon the bistability of those phase portraits. It should be noted that the same construction can be repeated on the restorative phase portraits in Fig.} \ref{FIG 3} {right. However those phase portraits are monostable and their ultra-slow modulation leads to a} {\it slow tonic spiking} ({\it i.e.} {a single spike necessarily followed by a rest period). This attractor differs from a bursting attractor by the absence of a bistable range in the bifurcation diagrams of Fig.} \ref{FIG 4}. {It can be shown that the persistence of (rest-spike) bistability in the singular limit is a hallmark of regenerative excitability} (Fig. \ref{FIG 3} left) {and that it cannot exist in restorative excitability} (Fig. \ref{FIG 3} right). See \cite{Franci2013} for a mode detailed discussion.
%{leading to a} {\it tonic spiking} {attractor. The absence of a regenerative saddle point lets any bistability of the restorative slow-fast subsystem vanish near the singular limit} \cite{Franci2013}. {Ultra-slow modulation of the unfolding parameter near the rest-spiking transition leads in this case to slow tonic spiking, in which each spike is necessarily followed by a rest period.
{ Modulation in} (\ref{EQ: 3D cusp dynamics}) {of the bifurcation parameter across the transcritical bifurcation of Fig.} \ref{FIG 3} {therefore provides a geometric transition from the slow tonic spiking attractor to the bursting attractor. This transition organizes the geometric route into bursting discussed in the next section}.

\section{{A physiological route to bursting}}
\label{SEC: the routes}

\subsection{A minimal three-time scale bursting model}

The recent paper \cite{Franci2012} introduces the planar neuron model
\begin{IEEEeqnarray}{rCl}\label{EQ: SIADS}
\dot V&=&V-\frac{V^3}{3}-n^2+I\IEEEyessubnumber\\
\dot n&=&\varepsilon(n_\infty(V-V_0)+n_0-n)\IEEEyessubnumber
\end{IEEEeqnarray}
Its phase portrait was shown to contain the pitchfork of Figure \ref{FIG 3} as an organizing center, leading to distinct types of excitability for distinct values of the unfolding parameters. The analysis of the previous section suggests that a bursting model is naturally obtained by augmenting the planar model (\ref{EQ: SIADS}) with ultra slow adaptation:
\begin{IEEEeqnarray}{rCl}\label{EQ: SIAM burst quantitative}
\dot V&=&kV-\frac{V^3}{3}-(n+n_0)^2+I-z\IEEEyessubnumber\\
\dot n&=&\varepsilon_n(V)\left(n_\infty(V-V_0) -n \right) \IEEEyessubnumber\\
\dot z&=&\varepsilon_z(V)(z_\infty(V-V_1)-z)\IEEEyessubnumber
\end{IEEEeqnarray}
Model (\ref{EQ: SIADS}) is essentially model (\ref{EQ: SIAM burst quantitative}) for $k=1$ and $z=0$, modulo a translation $n\gets n+n_0$. The dynamics (\ref{EQ: SIAM burst quantitative}b-\ref{EQ: SIAM burst quantitative}c) mimic the kinetics of gating variables in conductance-based models, where the steady-state characteristics $n_\infty(\cdot)$ and $z_\infty(\cdot)$ are monotone increasing (typically sigmoidal) and the time scaling $\varepsilon_n(\cdot)$ and $\varepsilon_z(\cdot)$ are Gaussian-like strictly positive functions.
{Details of model} (\ref{EQ: SIAM burst quantitative}) {for the numerical simulations of the paper are provided in Appendix} \ref{SEC: parameters}.

{The slow-fast subsystem} (\ref{EQ: SIAM burst quantitative}a\ref{EQ: SIAM burst quantitative}b) {shares the same geometric structure as} (\ref{EQ: cusp dynamics}). {After a translation $V\leftarrow V+V_0$, the right hand side of} (\ref{EQ: SIAM burst quantitative}a) {can easily be shown to be a universal unfolding of the winged cusp and the slow dynamics} (\ref{EQ: SIAM burst quantitative}b) {modulates its bifurcation parameter. Plugging the ultra-slow dynamics} (\ref{EQ: SIAM burst quantitative}c){, one recovers the same structure as} (\ref{EQ: 3D cusp dynamics}). {Therefore, the conclusions of Theorems} \ref{THM: cusp bist} and \ref{THM: cusp burst} {apply to} (\ref{EQ: SIAM burst quantitative}).

The difference between (\ref{EQ: SIAM burst quantitative}) and (\ref{EQ: 3D cusp dynamics}) is that the model (\ref{EQ: SIAM burst quantitative}) has the physiological interpretation of a reduced conductance-based model, with $V$ a fast variable that aggregates the membrane potential with all fast gating variables, $n$ a slow recovery variable that aggregates all the slow gating variables regulating neuronal excitability, and $z$ an ultra-slow adaptation variable that aggregates the ultra-slow gating variables that modulate the cellular rhythm over the course of many action potentials. Finally, $I$ models an external applied current.

\subsection{Model parameters and their physiological interpretation}
\label{SSEC: parameter interpretation}

\subsubsection*{The bifurcation parameter $n_0$ models the balance between restorative and regenerative ion channels}

\begin{figure}
\center
\includegraphics[width=0.9\textwidth]{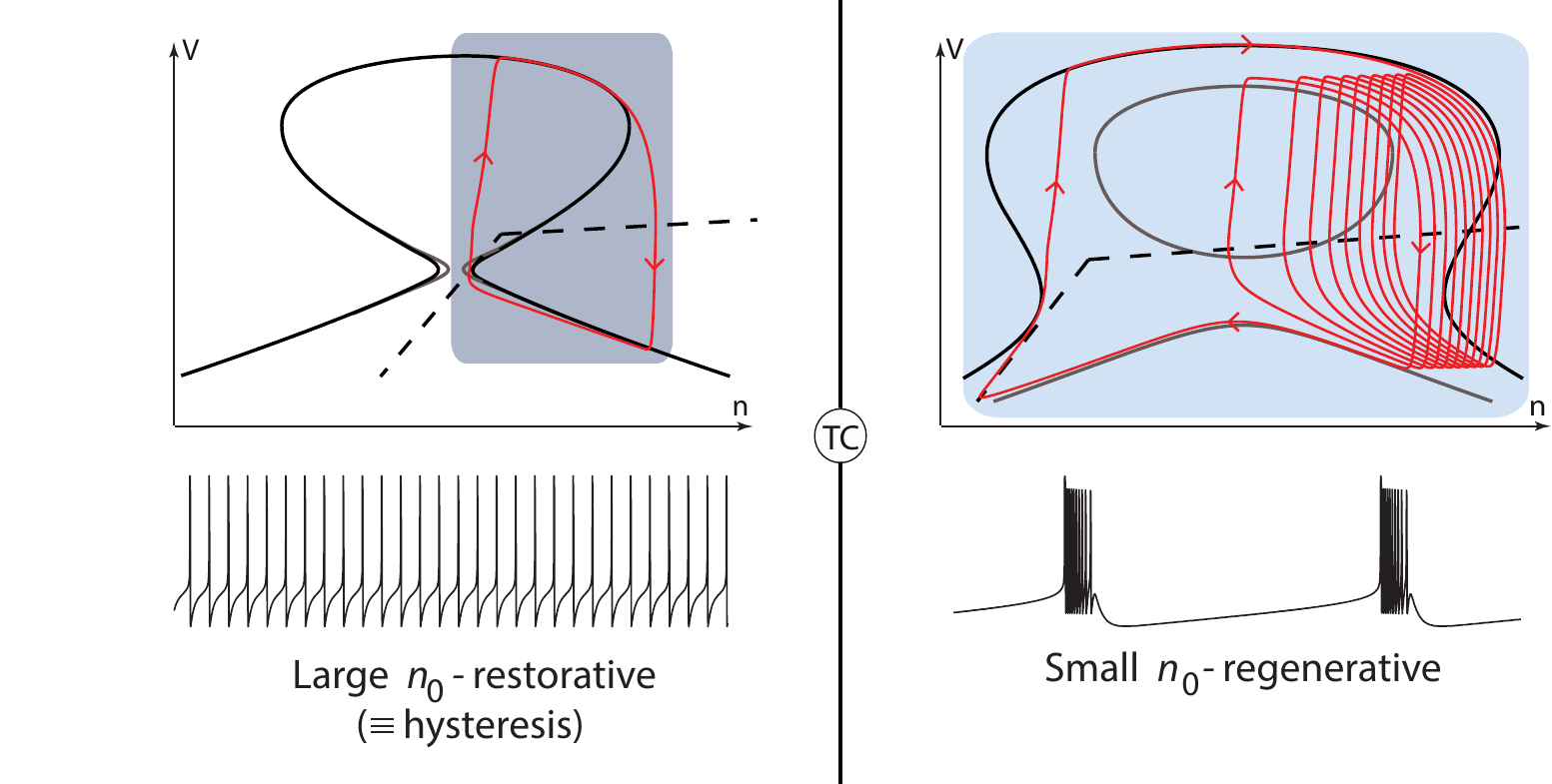}
\caption{{\bf Transition from restorative excitability (tonic firing) to regenerative excitability (bursting) in model (\ref{EQ: SIAM burst quantitative}) by sole variation of the bifurcation parameter $n_0$.} The analytical expression of the steady state functions $n_\infty(\cdot)$ and $z_\infty(\cdot)$ and numerical parameter values are provided in Appendix \ref{SEC: parameters}. The time scale is the same in the left and right {time series}.}\label{FIG 6}
\end{figure}

The central role of the bifurcation parameter $n_0$ in (\ref{EQ: SIAM burst quantitative}) was analyzed in \cite{Franci2012,Franci2013} and is illustrated in Fig. \ref{FIG 6}. The transcritical bifurcation variety in Fig. \ref{FIG 3} corresponds to the physiologically relevant transition from restorative excitability (large $n_0$) to regenerative excitability (small $n_0$). When the excitability is restorative, the recovery variable $n$ provides negative feedback on membrane potential variations near the resting equilibrium, a physiological situation well captured by FitzHugh-Nagumo model (or the hysteresis singularity). In contrast, when excitability is regenerative, the recovery variable $n$ provides positive feedback on membrane potential variations near the resting potential, a physiological situation that requires the quadratic term in (\ref{EQ: SIAM burst quantitative}a) (or the winged cusp singularity).

The value of $n_0$ in a conductance-based model reflects the balance between restorative and regenerative ion channels that regulate neuronal excitability. How to determine the balance in an arbitrary conductance-based model is discussed in \cite{Franci2013}. Note that the restorative or regenerative nature of a particular ion channel in the slow time-scale is an intrinsic property of the channel. A prominent example of restorative channel is the slow potassium activation shared by (almost) all spiking neurons. A prominent example of regenerative channel is the slow calcium activation encountered in most bursting neurons. The presence of regenerative channels in neuronal bursters is well established in neurophysiology. See e.g. \cite{Krahe2004,Astori2011}.

\subsubsection*{The affine unfolding parameter provides bursting by ultra-slow modulation of the current across the membrane}

For small $n_0$, the modulation of the ultra-slow variable $z$ creates a hyperbolic bursting attractor through the hysteretic loop described in Fig. \ref{FIG 4}. The burster becomes a single-spike limit cycle (tonic firing) for large $n_0$ (restorative excitability), that is, in the absence of rest-spike bistability in the planar model.

The presence of ultra-slow currents in neuronal bursters is well established in neurophysiology (see e.g. \cite{Astori2011}). A prominent example is provided by   ultra-slow calcium activated potassium channels.

\subsubsection*{Half activation potential affects the route to bursting}

The role of the unfolding parameter $\gamma$ in (\ref{EQ: 3D cusp dynamics}) is illustrated in Fig. \ref{FIG 3}: it provides two qualitatively distinct paths connecting the restorative and regenerative phase portraits. This role is played by the parameter $V_0$ in the planar model (\ref{EQ: SIADS}) studied in \cite{Franci2012}, which has the physiological interpretation of a half activation potential. The role of half-activation potentials in neuronal excitability is well documented in neurophysiology (see e.g. \cite{Putzier2009}). The role of this unfolding parameter in the route to bursting is discussed in the next subsection.

\subsubsection*{No spike without fast autocatalytic feedback}

The role of the unfolding parameter $k$ in (\ref{EQ: SIAM burst quantitative}) is to provide positive (autocatalytic) feedback in the fast dynamics. The prominent source of this feedback in conductance-based models is the fast sodium activation. It is well acknowledged in neurodynamics \cite{Izhikevich2007}.

The reduced model (\ref{EQ: SIAM burst quantitative}) makes clear predictions about its dynamical behavior in the absence of this feedback ({\it i.e.} $k=0$). Those predictions { are further discussed in Section} \ref{SSEC: burst across bursting types} and are in closed agreement with the experimental observation of ``small oscillatory potentials" when sodium channels are {shut down} with pharmacological blockers \cite{Guzman2009,Zhan1999} {or are poorly expressed during neuronal cell development} \cite{Liu1998}.

\subsection{A physiological route to bursting}
\label{SSEC: routes to burst}

\begin{figure}[h!]
\center
\includegraphics[width=\textwidth]{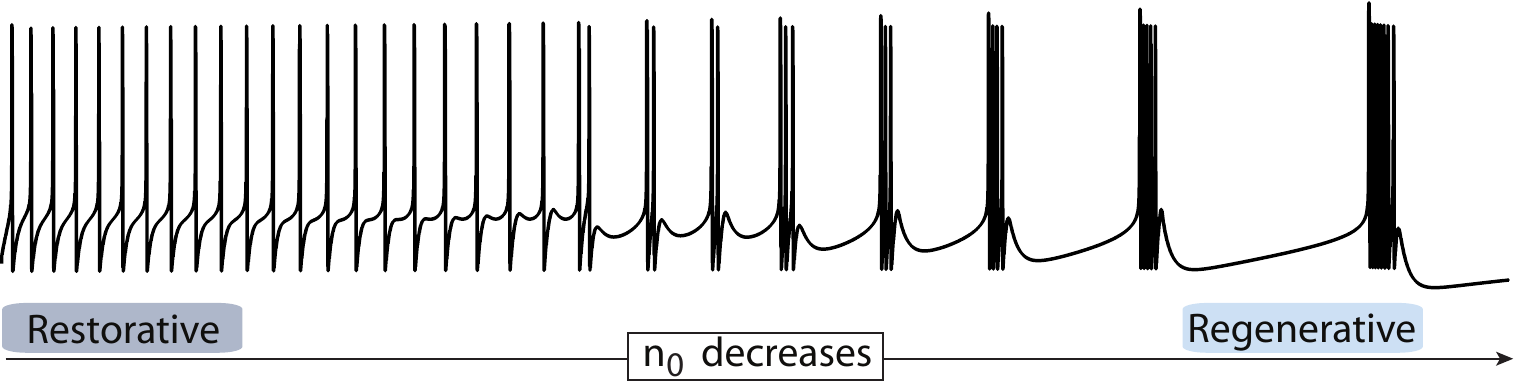}
\caption{{\bf Route from tonic firing to bursting in model (\ref{EQ: SIAM burst quantitative}) via a smooth variation of the bifurcation parameter $n_0$.} Rest of the parameter as in Figure \ref{FIG 6}.}\label{FIG 7}
\end{figure}

A central insight of the reduced model (\ref{EQ: SIAM burst quantitative}) is that it provides a route to bursting: fixing all unfolding parameters and varying only the bifurcation parameter $n_0$ leads to a smooth transition from tonic firing to bursting, see Fig. \ref{FIG 7}.
Smooth and reversible transitions between those two rhythms have been observed in many experimental recordings \cite{Sherman2001,Viemari2006}, making the route to burst an important signaling mechanism. The fact that the modulation is achieved simply through the bifurcation parameter $n_0$, {\it i.e.} the balance between restorative and regenerative channels, is of physiological importance because it is consistent with the physiology of experimental observations of routes into bursting \cite{Sherman2001,Viemari2006,Beurrier1999}.

The analysis in the above sections shows that the transition from single spike to bursting is through the transcritical bifurcation variety in model (\ref{EQ: cusp dynamics}). Looking at the singular limit $\varepsilon=0$ of (\ref{EQ: cusp dynamics}) near this transition variety provides further insight on the geometry of the route that leads to the appearance of the saddle-homoclinic bifurcation organizing the bistable phase-portrait. This route is organized by the path through the pitchfork bifurcation, which provides the most symmetric path across the transcritical variety. The generic transitions are understood by perturbing the degenerate path.

\begin{figure}[h!]
\center
\includegraphics[width=\textwidth]{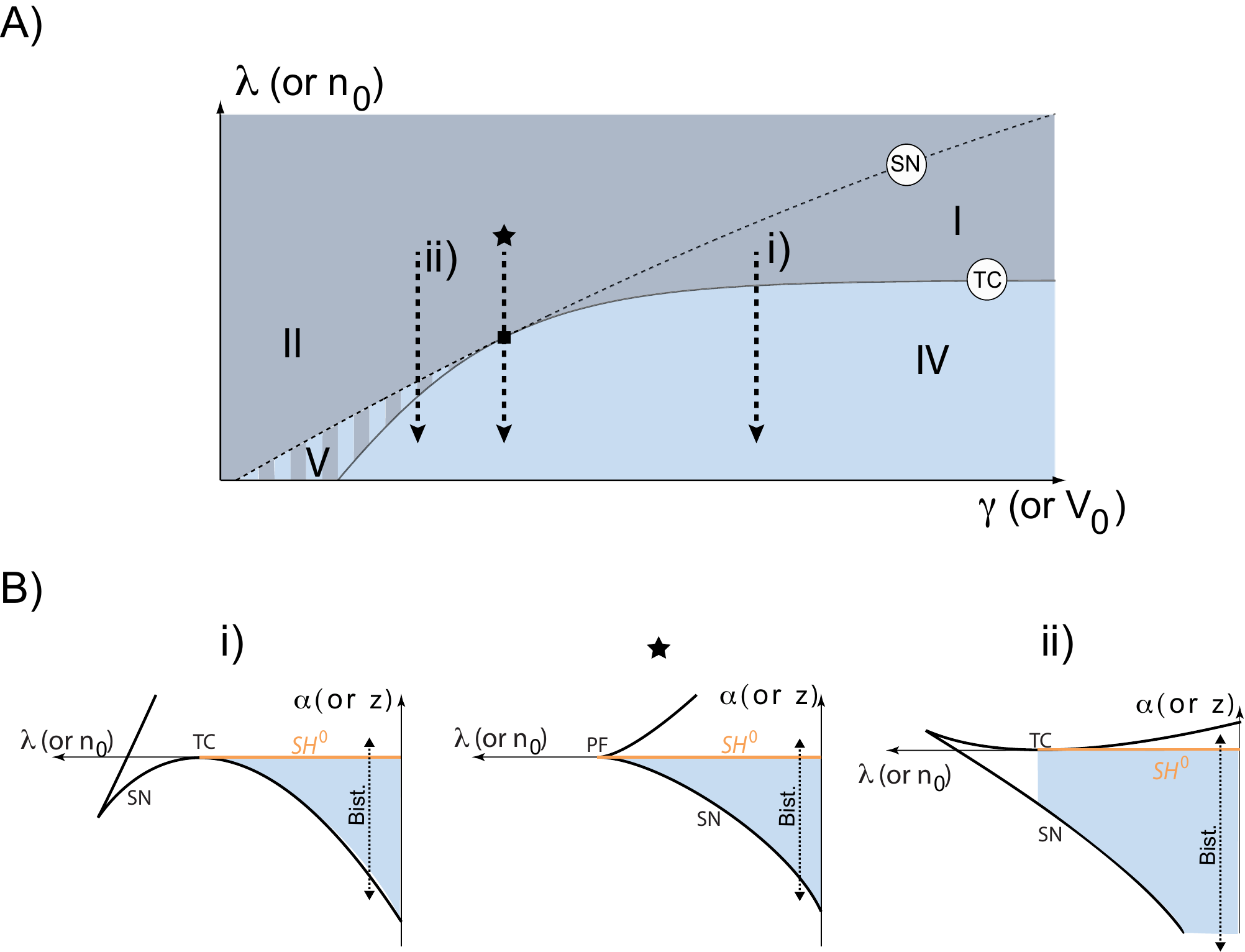}
\caption{{\bf Routes into bursting in the universal unfolding of the pithfork bifurcation.} {\bf A.} Qualitative projection of routes into bursting onto the $(V_0,n_0)$ (resp. $(\gamma,\lambda)$) of model (\ref{EQ: SIADS}) (resp. (\ref{EQ: cusp dynamics})) for $I=\frac{2}{3}$ (resp. $\alpha=\alpha_{TC}(\beta,\gamma)$, see Appendix \ref{SEC: trans variety}). Excitability is restorative in subregions $I$ and $II$, mixed in subregion $V$, and regenerative in subregion $IV$. See \cite{Franci2012} and \cite{Franci2013} for details concerning the underlying bifurcation mechanisms. The transition path labeled with a star depicts the degenerate path across the pitchfork. The generic paths {\it i)} and {\it ii)} are distinguished by different half activation potentials $V_0$ (resp. unfolding parameter $\gamma$). {\bf B.} Unfolding of transition paths in A along the $I$ (resp. $\alpha$) direction. Black thick lines denote branches of saddle-node (SN) bifurcation. In paths {\it i)} and {\it ii)}, the model undergoes a transcritical bifurcation (TC) as the path touches tangentially a branch of SN bifurcations. In the degenerate path, the model undergoes a pitchfork (PF) bifurcation as the path enters the cusp tangentially to both branches of SN bifurcations. The singular saddle-homoclinic loop, geometrically constructed in Figs. \ref{FIG 4} and \ref{FIG 9}, is denoted by $SH^0$ and determines the appearance of a singular bistable range persisting away from singular limit.
}\label{FIG 8}
\end{figure}

\begin{figure}[h!]
\center
\includegraphics[width=0.9\textwidth]{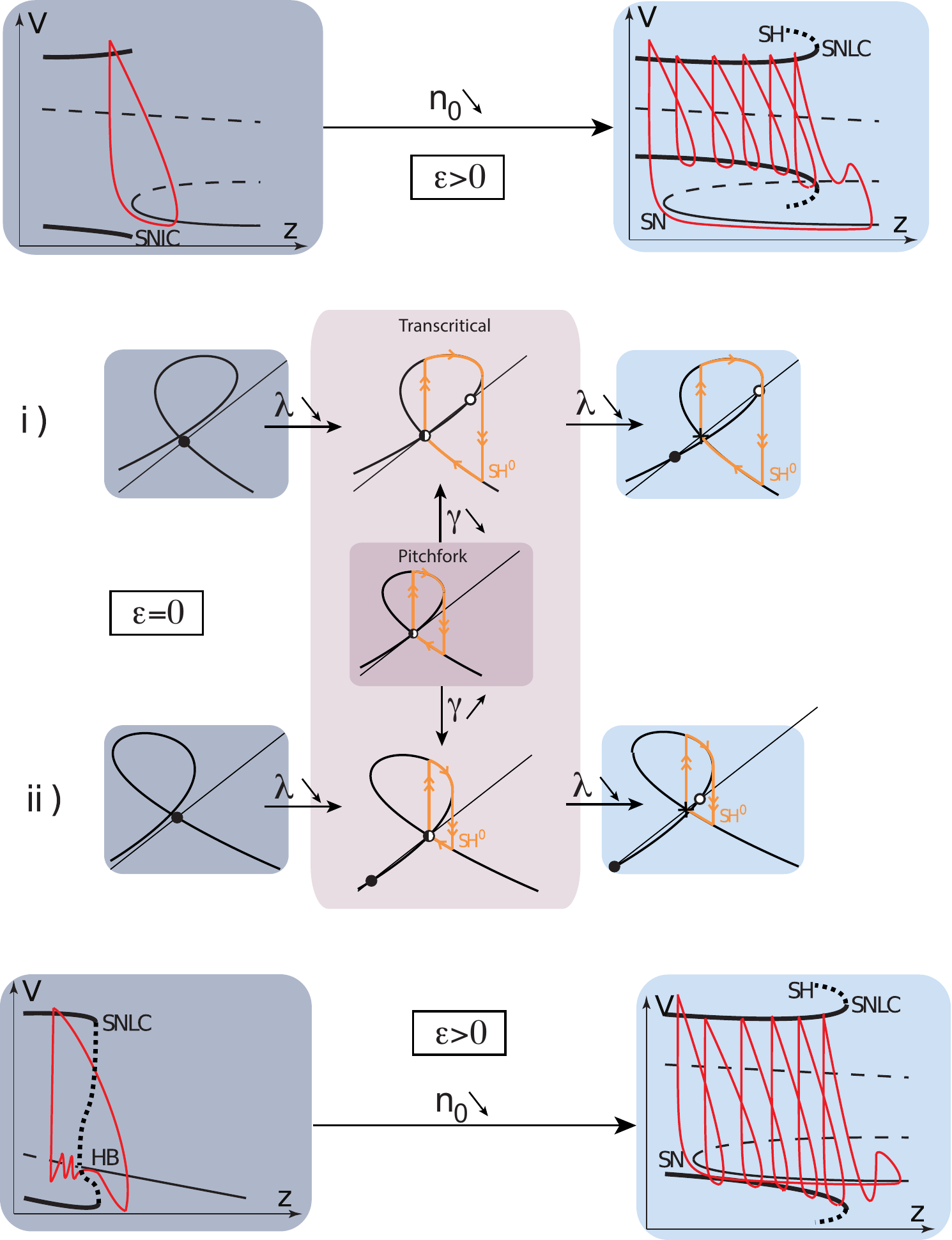}
\caption{{\bf Geometry of the two generic routes into bursting in the unfolding of the pitchfork bifurcation in model (\ref{EQ: 3D cusp dynamics}) and model (\ref{EQ: SIAM burst quantitative}).}}\label{FIG 9}
\end{figure}

Fig. \ref{FIG 8}A shows the qualitative projection of those paths onto the $(V_0,n_0)$ parameter chart obtained in model (\ref{EQ: SIADS}) for $I=\frac{2}{3}$. The chart is reproduced from \cite{Franci2012}. The same qualitative picture is obtained for the $(\gamma,\lambda)$ parameter chart of the abstract model (\ref{EQ: cusp dynamics}) at $\alpha=\alpha_{TC}(\beta,\gamma)$ (see Appendix \ref{SEC: trans variety}). The chart associates different excitability types (as well as their restorative or regenerative nature, see \cite{Franci2013}) to distinct bifurcation mechanisms. Unfolding those paths along the $I$ (or $\alpha$) direction leads to the bifurcation diagrams in Fig. \ref{FIG 9}B. They reveal (in the singular limit) the onset of the bistable range organized by the singular saddle-homoclinic loop $SH^0$ as paths cross the transcritical bifurcation variety.
%They show the appearance of the singular saddle-homoclinic loop $SH^0$ (geometrically constructed in Figs. \ref{FIG 4} and \ref{FIG 9}), as the paths cross the transcritical bifurcation variety, and the associated extended (singular) bistable range. In \cite{FRDRSE_balance}, different conductance based models were shown to travel the same qualitative paths via the variation of physiological parameters.

The same qualitative picture persists for $\varepsilon>0$.  Fig. \ref{FIG 9} illustrates how the appearance of the singular saddle-homoclinic loop is accompanied, for $\varepsilon>0$, by a smooth transition from a monostable (SNIC - route {\it i)} ) or barely bistable (sub. Hopf - route {\it ii)} ) bifurcation diagram to the robustly bistable bifurcation diagram constructed in the sections above (Fig. \ref{FIG 5}). Through ultra-slow modulation of the unfolding parameter $\alpha$, this transition geometrically captures the transition from tonic spiking to bursting via the sole variation of the bifurcation parameter.

{The strong agreement between the mathematical insight provided by singularity theory and the known electrophysiology of bursting is a peculiar feature of the proposed approach. There is a direct correspondence between the bifurcation and unfolding parameters of the winged cusp and the physiological minimal ingredients of a neuronal burster. In particular, our analysis predicts that any bursting neuron must possess at least one physiologically regulated slow regenerative channel. This prediction needs to be tested systematically but we have found no counter-example in the bursting neurons we have analyzed to date.}

\section{{Normal form reduction of conductance-based models}}
\label{SEC: w cusp in cb modes}

\subsection{A two dimensional reduction}
\label{SSEC: cb models two d red}

The winged cusp singularity emerges as an organizing center of rhythmicity  in the reduced neuronal model (\ref{EQ: SIAM burst quantitative}), but a legitimate question is whether this singularity can be traced in arbitrary (high-dimensional) conductance-based models. Our recent paper \cite{Franci2013} addresses a closely related question for the transcritical variety. It provides an analog of the bifurcation parameter $n_0$ in arbitrary conductance-based models of the form
\begin{IEEEeqnarray}{rCll}\label{EQ: generic cb model}
C_m\dot V&=&-\sum_{\iota}\bar g_\iota m_\iota^{a_\iota} h_\iota^{b_\iota} (V-E_\iota)+I_{app},\nonumber\\
&=:&I_{ion}(V,x^f,x^s,x^{us})+I_{app}&\IEEEyessubnumber\\
\tau_{x^{f}_j}(V)\dot x^f_j&=&-x^f_j+x^f_{j,\infty}(V),&\quad j=1,\ldots,n_f\IEEEyessubnumber\\
\tau_{x^s_j}(V)\dot x^s_j&=&-x^s_j+x^s_{j,\infty}(V),&\quad j=1,\ldots,n_s\IEEEyessubnumber\\
\tau_{x^{us}_j}(V)\dot x^{us}_j&=&-x^{us}_j+x^{us}_{j,\infty}(V),&\quad j=1,\ldots,n_{us}\IEEEyessubnumber
\end{IEEEeqnarray}
where $\iota$ {runs through} all ionic currents, $x^f:=[x^f_j]_{j=1,\ldots,n_f}$ denotes the $n_f$-dimensional column vector of fast gating variables, $x^s:=[x^s_j]_{j=1,\ldots,n_s}$ denotes the $n_s$-dimensional column vector of slow gating variables, and $x^{us}:=[x^{us}_j]_{j=1,\ldots,n_{us}}$ denotes the $n_{us}$-dimensional column vector of ultra-slow variables (see also \cite{Franci2013} for more details on the adopted notation).

Following common analysis methods in neurodynamics, we want to reduce the (possibly) high-dimensional model (\ref{EQ: generic cb model}) to a two-dimensional model of the form
\begin{IEEEeqnarray}{rCl}\label{EQ: CB model 2d red}
\dot V
%&=&I_{ion}\left( V,x_\infty^f(V), \tilde x^s(n),
%\big|_{\begin{subarray}{l}x^s_1=n,\\
%x^s_i=x^s_{i,\infty}\left({n_{,\infty}^s}^{-1}(n)\right),\\ i= 2,\ldots,n_s\end{subarray}},
%\bar x^{us}\right)+I_{app}\nonumber\\
&=&F(V,n)+I\IEEEyessubnumber\\
\tau(V)\dot n^s&=&-n+n_\infty^s(V)\IEEEyessubnumber
\end{IEEEeqnarray}
where $V$ is the fast voltage and $n$ is a slow aggregate variable. We achieve this reduction by first considering the singular limit of three time scales leading to a quasi-steady state approximation for fast gating variables, that is
\begin{equation}\label{EQ: cb inst fast}
x^f_j\equiv x^f_{j,\infty}(V),
\end{equation}
for all $ j=1,\ldots,n_f$, and freezing ultra-slow variables, that is setting
\begin{equation*}%\label{EQ: cb freeze ultra-slow}
x^{us}_j\equiv \bar x^{us}_j,
\end{equation*}
for all $j=1,\ldots,n_{us}$, where the values $\bar x^{us}_j$ belong to the physiological range of the different variables. The remaining dynamics read as
\begin{IEEEeqnarray*}{rCl}%\label{EQ: slow fast cb}
\dot V&=&I_{ion}(V,x^f_\infty(V),x^s,\bar x^{us})+I_{app},\\
\tau(V)\dot x^s_j&=&(-x^s_j+x^s_{j,\infty}(V)),\quad j=1,\ldots,n_s
\end{IEEEeqnarray*}
which is a fast-slow system with $V$ as fast variable and $x^s$ as slow variables.

The planar reduction proceeds from the change of variables
\begin{IEEEeqnarray*}{rCl}%\label{EQ: cb change of variables}
n&=&x_1^s,\\
n_i^\perp&=&x_i^s-x_{i,\infty}^s(n_\infty^{-1}(n)),\quad i=2,\ldots,n_s,
\end{IEEEeqnarray*}
This change of variable is globally invertible by monotonicity of the (in)activation functions $x_{i,\infty}^s$. Under the additional simplifying assumption of identical time constants
\begin{equation}\label{EQ: cb homo time slow}
\tau_{x^s_j}(V)= \tau(V)\geq\epsilon^{-1}\gg 1,
\end{equation}
for all $V\in\mathbb R$ and all $j=1,\ldots,n_s$, it is an easy calculation to show that
\begin{IEEEeqnarray*}{rCl}%\label{EQ: n perp exp converg}
\tau(V^\star)\dot n_i^\perp&=&- n_i^\perp +\mathcal O((n-n^\star)^2,(n-n^\star)(V-V^\star),(V-V^\star)^2)
\end{IEEEeqnarray*}
around any equilibrium $(V^\star,(x_i^s)^\star):=(V^\star,x^s_{i,\infty}(V^\star))$. It follows that, locally around {\it any} equilibrium, the two dimensional manifold
\begin{IEEEeqnarray*}{rCl}
\mathcal M_{red}&:=&\left\{(V,x^s_\infty)\in\mathbb R\times[0,1]^{n_s}:\ n_i^\perp=0,\ i=2,\ldots,n_s \right\}\\
&=&\left\{(V,x^s_\infty)\in\mathbb R\times[0,1]^{n_s}:\ x_i^s=x_{i,\infty}^s(n_\infty^{-1}(n)),\ i=2,\ldots,n_s \right\}
\end{IEEEeqnarray*}
is exponentially attractive.

It should be stressed that the (harsh) simplifying assumption (\ref{EQ: cb homo time slow}) is necessary only around the steady-state value $V^\star$ and that the hyperbolic decomposition is robust to small perturbations \cite{Hirsch1977}. It should also be observed that the proposed two-dimensional reduction is a straightforward generalization of the classical two-dimensional reduction of Hodgkin-Huxley model \cite{FitzHugh1961,Rinzel1985} that rests on setting sodium activation to steady state ($m_{Na}\equiv m_{Na,\infty}(V)$) and using an algebraic relationship between the sodium inactivation and the potassium activation (usually in the form $h\simeq1-n$).

\subsection{The winged cusp planar model (\ref{EQ: cusp dynamics}) is a local normal form of slow-fast conductance based models}
\label{SSEC: wcusp in cb models}

Given an equilibrium $(V^\star,n_\infty(V^\star))$ of (\ref{EQ: CB model 2d red}), consider the  (linear) change of variables
\begin{IEEEeqnarray*}{rCl}\label{EQ: CB model 2d red varia variables}
x&=&V-V^\star\\
y&=&\frac{n-n_\infty(V^\star)}{\frac{\partial n_\infty}{\partial V}(V^\star)}
\end{IEEEeqnarray*}
The $y$ dynamics is particularly simple. Indeed, by simple Taylor expansion,
\begin{IEEEeqnarray*}{rCl}%\label{EQ: CB model 2d red y dyna}
\dot y&=&
%\left(\frac{\partial n_\infty}{\partial V}(V^\star) \right)^{-1}\frac{1}{\tau(V^\star)}(-n+n_\infty(V))\nonumber\\
%&=&\varepsilon\left(\frac{-n+n_\infty(V^\star)}{\frac{\partial n_\infty}{\partial V}(V^\star)}+\frac{\frac{\partial n_\infty}{\partial V}(V^\star) }{\frac{\partial n_\infty}{\partial V}(V^\star) }(V-V^\star) \right)+\mathcal O((V-V^\star)^2)\nonumber\\
%&=&
\varepsilon(x-y)+\mathcal O(x^2),\\
\varepsilon&:=&\frac{1}{\tau(V^\star)}\ll 1.
\end{IEEEeqnarray*}
In the new coordinates, (\ref{EQ: CB model 2d red}) reads
\begin{IEEEeqnarray}{rCl}\label{EQ: CB model 2d red varia dyn}
\dot x&=&F\left(x+V^\star,n_\infty(V^\star)+\frac{\partial n_\infty}{\partial V}(V^\star)y\right)+I\IEEEyessubnumber\\
\dot y&=&\varepsilon(x-y)+\mathcal O(x^2)\IEEEyessubnumber.
\end{IEEEeqnarray}
Simple computations show that (\ref{EQ: CB model 2d red varia dyn}a) satisfies
\begin{IEEEeqnarray*}{rCl}%\label{EQ: CB model 2d red dx expressions}
\frac{\partial\dot x}{\partial x}(0,0)&=&\frac{\partial I_{ion}}{\partial V}+\sum_{i=1}^{n_f}\frac{\partial I_{ion}}{\partial x_i^f}\frac{\partial x^f_{i,\infty}}{\partial V}\\
\frac{\partial\dot x}{\partial y}(0,0)&=&\sum_{i=1}^{n_s}\frac{\partial I_{ion}}{\partial x_i^s}\frac{\partial x^s_{i,\infty}}{\partial V}
\end{IEEEeqnarray*}
where the right hand sides are intended computed at $V=V^\star,x^f=x^f_\infty(V^\star),x^s=x^s_\infty(V^\star)$, and $x^{us}=\bar x^{us}$.

We claim that the critical manifold $\dot x=0$ of (\ref{EQ: CB model 2d red}) has a degenerate singularity provided that
\begin{itemize}
\item[] {\it (i)} the full slow-fast subsystem has a degenerate equilibrium, that is, the Jacobian of the slow-fast subsystem (\ref{EQ: generic cb model}a-\ref{EQ: generic cb model}c) is singular
\item[] {\it (ii)} at such equilibrium, the contributions of slow restorative and slow regenerative channels \cite{Franci2013} are perfectly balanced, that is
\begin{IEEEeqnarray*}{rCl}%\label{EQ: generic balance}
\sum_{i=1}^{n_s}\frac{\partial I_{ion}}{\partial x_i^s}\frac{\partial x^s_{i,\infty}}{\partial V}&=&0.
\end{IEEEeqnarray*}
\end{itemize}
To prove our claim we notice with similar computations as \cite{Franci2013} that conditions {\it (i)} and {\it (ii)} imply that $$\frac{\partial I_{ion}}{\partial V}+\sum_{i=1}^{n_f}\frac{\partial I_{ion}}{\partial x_i^f}\frac{\partial x^f_{i,\infty}}{\partial V}=0,$$ which is equivalent to the Jacobian of the fast subsystems (\ref{EQ: generic cb model}a-\ref{EQ: generic cb model}b) being singular. Hence, when conditions {\it (i)} and {\it (ii)} are fulfilled,
\begin{equation}\label{EQ: varia dyn deg cond}
\frac{\partial\dot x}{\partial x}(0,0)=\frac{\partial\dot x}{\partial y}(0,0)=0.
\end{equation}

Property (\ref{EQ: varia dyn deg cond}) ensures that the critical manifold of (\ref{EQ: CB model 2d red varia dyn}) has a codimension$>0$ singularity at the origin (where, as usual, the slow variable $y$ plays the role of the bifurcation parameter). This singularity corresponds to the transcritical bifurcation detected in arbitrary conductance based models in \cite{Franci2013}. It is indeed proved in \cite{Franci2013} that conditions {\it (i)} and {\it (ii)} enforce a transcritical bifurcation in the associated conductance based model.

Algebraically, (\ref{EQ: varia dyn deg cond}) ensures that, similarly to the bifurcation parameter in the winged cusp universal unfolding (see Section \ref{SSEC: rest-spike in wcusp}), $y$ modulates non-monotonically the fast $x$ dynamics. Physiologically, it captures in the reduced model the non-monotone modulation of membrane potential dynamics by slow restorative (providing negative feedback) and slow regenerative (providing positive feedback) ion channels.

\begin{figure}[h!]
\center
\includegraphics[width=0.9\textwidth]{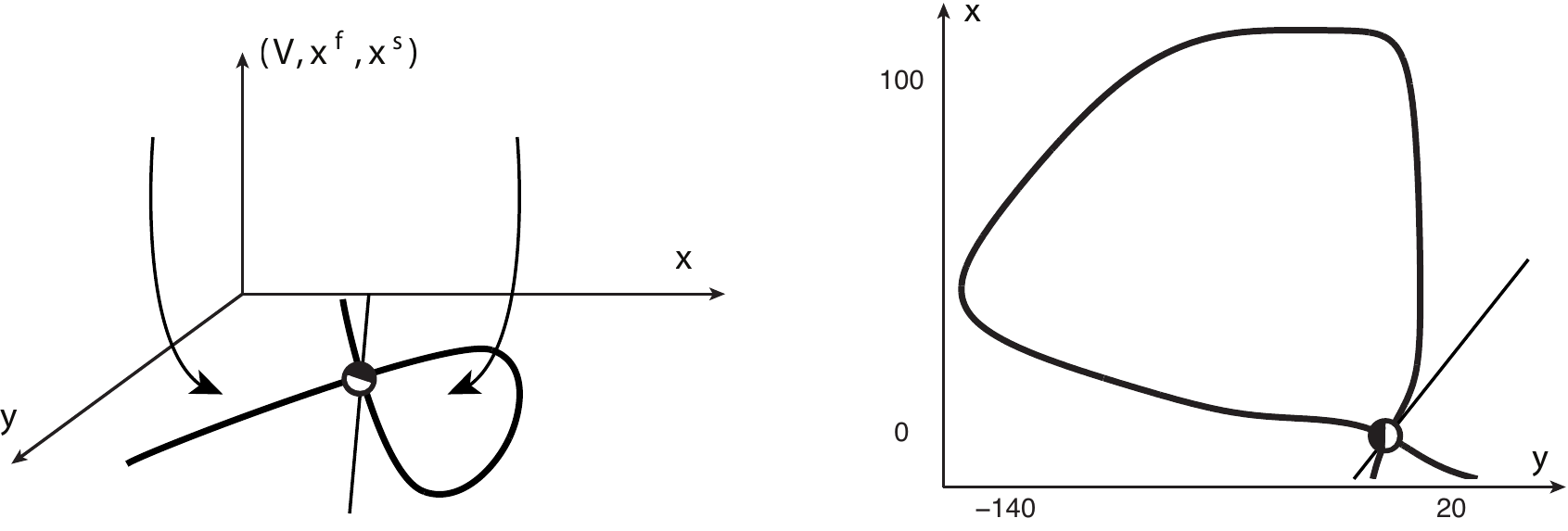}
\caption{A transcritical bifurcation in the universal unfolding of the winged cusp organizes the dynamics of the two-dimensional reduction of generic conductace based models. Left: Sketch of the dynamics on the two-dimensional invariant manifold $\mathcal M_{red}$. Right: construction of the two dimensional reduction (\ref{EQ: CB model 2d red varia dyn}) at the transcritical bifurcation in the Hodgkin-Huxley model augmented with a calcium current (\ref{EQ: HH model})-(\ref{EQ: HH CA dynamics}).}\label{FIG XX1}
\end{figure}

We use the algorithm in \cite{Franci2013} to detect the degenerate dynamics of (\ref{EQ: CB model 2d red varia dyn}) in arbitrary conductance based models. This construction reveals that the transcritical bifurcation is part of the transcritical transition variety in the universal unfolding of the winged cusp. The result is sketched in Figure \ref{FIG XX1} left and verified numerically in the Hodgkin-Huxley model augmented with a calcium current  in Figure \ref{FIG XX1} right. The model and its reduction are presented and further discussed in Section \ref{SSSE: winged-cusp phase plane in HHCa} below. The obtained phase plane is organized by the mirrored hysteresis bifurcation diagram of the normal form (\ref{EQ: cusp dynamics}) in Fig. \ref{FIG 2}, in the limiting case in which the two hystresis branches merge at the transcritical bifurcation. This provides an indirect proof that the global phase plane is organized by the winged cusp. This singularity is indeed the only (codimension$\leq3$) singularity exhibiting the mirrored hysteresis in its universal unfolding (see \cite[Section IV.4]{Golubitsky1985}).

One can push forward the singularity analysis and derive an algorithm to enforce the degenerate conditions of the winged cusp rather than the transcritical bifurcation by using additional model parameters as auxiliary parameters \cite[Section III.4]{Golubitsky1985}. This would lead to the conclusion that the critical manifold of the reduced dynamics (\ref{EQ: CB model 2d red varia dyn}) is actually a versal unfolding of the winged cusp. Alternatively, one can modulate model parameters and show that their variations recover all persistent bifurcation diagrams of the winged cusp. Such computations are however lengthy and bring no new information to the picture presented here.
%Below we simply pick some auxiliary parameters in the Hodgkin-Huxley model with a calcium current and show that their modulation uncovers the persistent bifurcation diagram of the winged cusp in the critical manifold of the associated reduced dynamics.

\subsection{Application to the Hodgkin-Huxley model augmented with a regenerative channel}
\label{SSSE: winged-cusp phase plane in HHCa}

The first conductance-based model appears in the seminal paper of Hodgkin-Huxley \cite{Hodgkin1952}
\begin{IEEEeqnarray}{rCl}\label{EQ: HH model}
C \dot V &=& - \bar g_Kn^4(V-V_K) - \bar g_{Na}m^3h(V-V_{Na})- g_l(V-V_l)+I\IEEEyessubnumber\\
\tau_m(V)\dot m&=&-m+m_\infty(V)\IEEEyessubnumber\\
\tau_n(V)\dot n&=&-n+n_\infty(V)\IEEEyessubnumber \\
\tau_h(V)\dot h&=&-h+h_\infty(V),\IEEEyessubnumber
\end{IEEEeqnarray}
where the time constants $\tau_x$ and the steady state characteristics $x_\infty$, $x=m,n,h$ are chosen in accordance with the original model (see Appendix \ref{SEC: parameters HH}). The model only accounts for two ionic currents: sodium, with its fast activation variable $m$ and slow inactivation $h$, and potassium, with slow activation $n$. The classical phase portrait reduction \cite{FitzHugh1961,Rinzel1985} is obtained with the quasi-steady state approximation $m\simeq m_\infty(V)$ and the empirical fit $h\simeq 1-n$. It is well known that in its physiological part ($0<n<1$) this phase portrait is qualitatively the FitzHugh phase portrait in Fig. \ref{FIG 1}. But we showed in \cite[Figure 5]{Drion2012} that the entire phase portrait ($n\in\mathbb R$) indeed also contains the ``mirrored" phase portrait of Fig. \ref{FIG 2}. This observation suggests that a winged cusp organizes the fast subsytem (\ref{EQ: HH model}a-\ref{EQ: HH model}b) of Hodkgin-Huxley dynamics. The singularity is found in a non-physiological range of the phase space ($n<0$), which is consistent with the absence of slow regenerative currents in the model.

The missing element in Hodgkin-Huxley model to make the winged cusp physiological is a slow regenerative ion channel. Following \cite{Drion2012}, we add the calcium current
\begin{IEEEeqnarray}{rCl}\label{EQ: HH CA dynamics}
I_{Ca,L}&=&-\bar g_{Ca}d(V-V_{Ca}) \IEEEyessubnumber\\
\tau_d(V) \dot d&=&-d + d_\infty(V). \IEEEyessubnumber
\end{IEEEeqnarray}
The algorithm in \cite{Franci2013} detects a transcritical bifurcation for
$$V^\star\simeq-61.2730,\quad g_{Ca}^\star\simeq0.2520,\quad I^\star\simeq-30.7694. $$
Following the construction in Section \ref{SSEC: wcusp in cb models}, in particular, Eq. (\ref{EQ: CB model 2d red varia dyn}), the associated reduced variational dynamics at the transcritical bifurcation reads
\begin{IEEEeqnarray*}{rCl}%\label{EQ: HH model normal form}
\dot x &=& - \bar g_K\left(n_\infty(V^\star)+y\frac{\partial n_\infty}{\partial V}(V^\star)\right)^4(V^\star+x-V_K)\nonumber\\
&&- \bar g_{Na}m_\infty(V^\star+x)^3\left(h_\infty(V^\star)+y\frac{\partial h_\infty}{\partial V}(V^\star)+\mathcal O(y^2)\right)(V-V_{Na})\\
&&-\bar g_{Ca}^\star\left(d_\infty(V^\star)+y\frac{\partial d_\infty}{\partial V}(V^\star)+\mathcal O(y^2)\right)(V-V_{Ca})\nonumber\\
&&- g_l(V-V_l)+I^\star\\
\dot y&=&\varepsilon(x-y)+\mathcal O(x^2).
\end{IEEEeqnarray*}
Its phase plane is drawn in Figure \ref{FIG XX1} right.

%Playing with the model parameters we can uncover in this phase plane the persistent bifurcation diagrams of the winged cusp (Fig. {\color{red} TO BE DRAWN}), which provides a numerical evidence that the model is a versal (since the number of auxiliary parameters exceeds the codimension of the singularity) unfolding of the winged cusp and further supports the role of the normal form (\ref{EQ: cusp dynamics}) in organizing the slow-fast dynamics of conductance-based models.

\begin{figure}[h!]
\center
\includegraphics[width=0.8\textwidth]{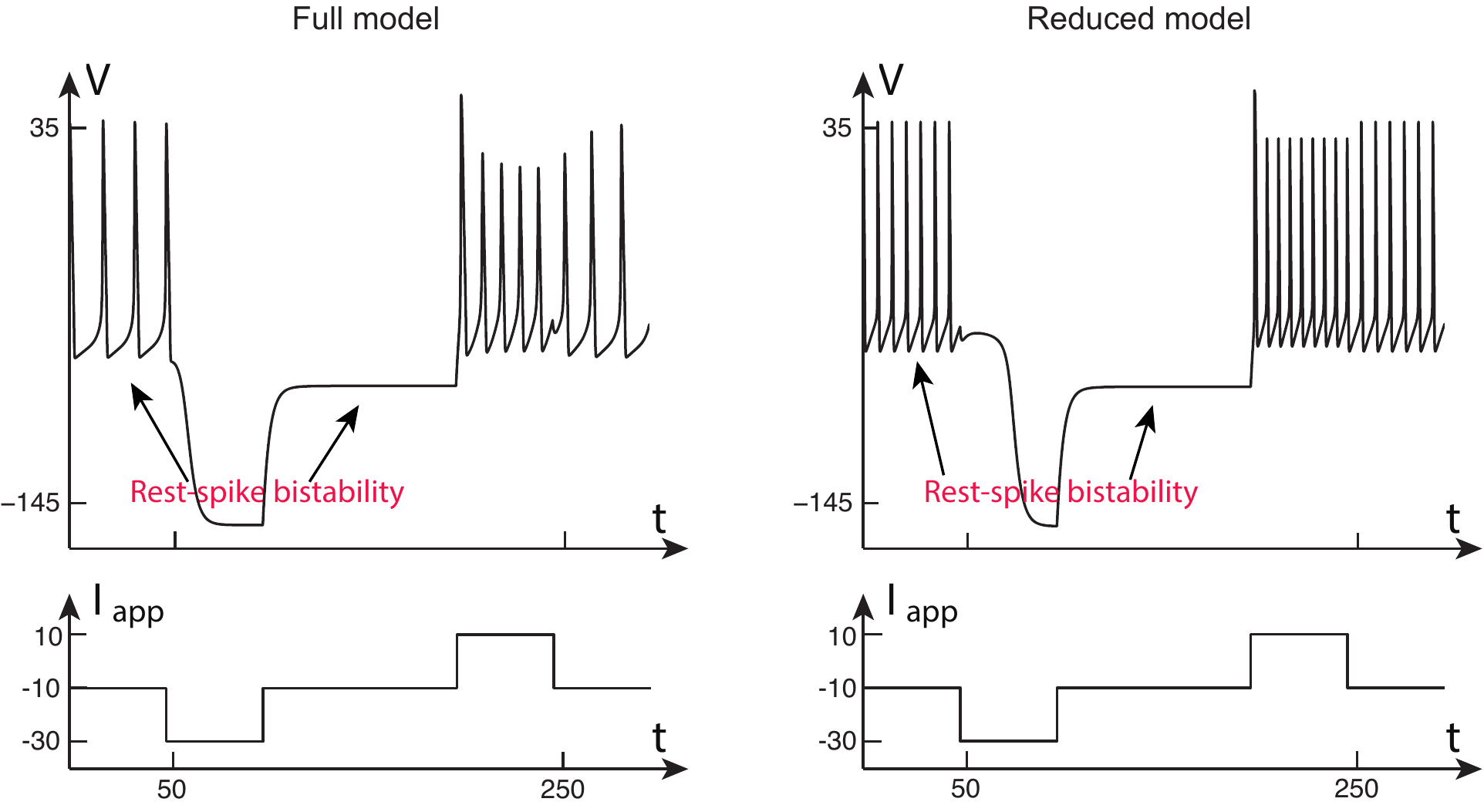}
\caption{Comparison of the full Hodgkin-Huxley augmented with a calcium current (\ref{EQ: HH model})-(\ref{EQ: HH CA dynamics}) and its two-dimensional reduction, obtained by applying the reduction procedure of Section \ref{SSEC: cb models two d red}.}\label{FIG XX2}
\end{figure}

We now apply the global two-dimensional reduction described in Section \ref{SSEC: cb models two d red}, in particular, Eq. (\ref{EQ: CB model 2d red}), to model (\ref{EQ: HH model}-\ref{EQ: HH CA dynamics}). To this aim, we express all variables in terms of potassium activation $n$. Since in the original model its activation function cannot be explicitly inverted, we use the exponential fitting
$$n_\infty(V)=\frac{1}{1+e^{0.06(11.6-V)}},\quad n_\infty^{-1}(n)=11.6-\frac{1}{0.06}\ln\left(\frac{1}{n}-1\right)$$
Figure \ref{FIG XX2} provides a comparison of the behavior of the original and reduced models. Despite quantitative differences (in particular, as in the reduction of the original Hodgkin-Huxley model, treating fast variables as instantaneous increases spiking frequency), the reduced model faithfully captures the qualitative behavior of its high-dimensional counterpart, for instance, rest-spike bistability. Phase plane analysis of the associated normal form (\ref{EQ: cusp dynamics}) provides a clear geometrical interpretation of such dynamical behavior (Fig. \ref{FIG 2}).

\subsection{The role of ultra-slow variables}
\label{SSEC: ultra slow unfolding}

Ultra-slow variables appear in a variety of forms: ultra-slow gating variables (e.g. inactivation of calcium channels), intracellular calcium (e.g. SK channels),  metabotropic regulation of channel expression (e.g. regulation of calcium channel expression by serotonin receptors), homeostatic regulation of channel expression (e.g. calcium dependent expression of ion channels), etc.~. As such, they do not allow a systematic analysis as for slow-gating variables. However, their effect on the model reduction (\ref{EQ: CB model 2d red}) can be understood in terms of modulation of the unfolding parameters of the associated normal form. The observation that the many (auxiliary) parameters of conductance based models might naturally provide a versal unfolding of the winged cusp organizing their fast critical manifold suggests that variations in ultra-slow variables act as ultra-slow modulation of the unfolding parameters in the associated normal form. The effect of ultra-slow variables is thus constrained to reshape the geometry of the slow fast phase portrait. This might lead to ultra-slow adaptation mechanisms (similarly to the action of $\alpha$ in Fig. \ref{FIG 4}) or to even slower modulation mechanisms (similarly to the action of $k$ and $V_0$ in Fig. \ref{FIG 11} below).

Clearly, this does not permit to conclude precise results on the global dynamics of a multi-timescale model, but suggest that the low dimensional bursting modulation mechanism described here has a strong relevance for generic conductance-based models.

\section{{Modulation of bursting by unfolding parameters and its physiological interpretation}}
\label{SEC: why the winged cusp}

\subsection{Bursting modeling and unfolding theory}

The rich literature on mathematical modeling of bursting calls for a few comparisons with the model proposed in the present paper. The geometry of our bursting attractor is the most classical one of a saddle-homoclinic burster (one out of the 16 bursting attractors in the recent classification of Izhikevich, see \cite[page 376]{Izhikevich2007}). Such an attractor is for instance found in the early bursting model of  Hindmarsh and Rose \cite{Hindmarsh1984}. The two models exhibit an analog geometry: the mirror of the classical Fitz-Hugh phase portrait, obtained here by mirroring the fast variable cubic nullcline, is obtained there by mirroring the monotone activation function of the recovery variable. But the Hindmarsh-Rose model lacks the organization of some high-codimension singularity, making it impractical for modulation studies  (see, e.g., \cite{Shilnikov2008})  and for physiological interpretability.

The more recent literature on bursting has certainly exploited unfolding theory around high-codimension bifurcations  to construct different types of bursting attractors. A non exhaustive list is \cite{Rinzel1987b,Bertram1995,Golubitsky2001,Izhikevich2000} and the references discussed in \cite[page 376]{Izhikevich2007}. The outcome of those studies is a useful mathematical classification between different bursting attractors organized by different bifurcations but it is not clear how to use this classification for modulation studies.  A possible reason is that most of those references construct bursting models from  restorative phase portraits that retain the qualitative organization of Fitz-Hugh model by a hysteresis singularity. Such models lack the transcritical bifurcation that organizes the normal form reduction of general bursting conductance based models.

The approach of the present paper differs from earlier studies in starting from the cusp singularity, inspired by our original observation that the mirrored hysteresis phase portrait organizes the reduced Hodgkin-Huxley dynamics \cite[Fig. 5]{Drion2012}. The direct link between the mathematical unfolding of the cusp singularity and the local normal form of conductance-based models in the vicinity of their transcritical bifurcation is probably crucial in using unfolding theory to understand the modulation of bursting in neuronal models.

\subsection{A geometrical and physiological modulation of a burster across bursting types}
\label{SSEC: burst across bursting types}

\begin{figure}[h!]
\center
\includegraphics[width=\textwidth]{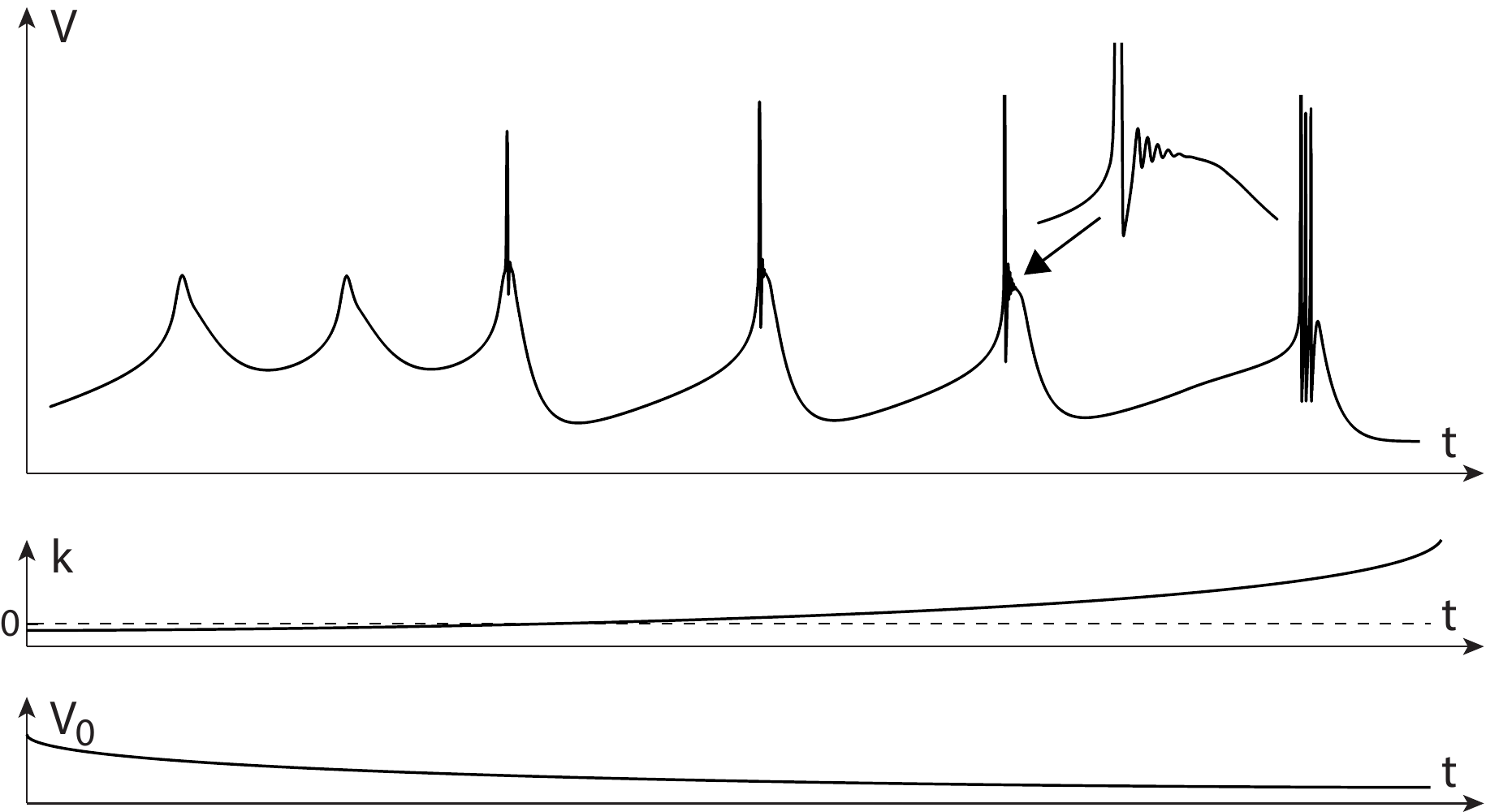}
\caption{{\bf Modulation of model (\ref{EQ: SIAM burst quantitative}) across different bursting wave forms.} Increasing the fast positive feedback gain $k$ and decreasing the half-activation potential $V_0$ the model smoothly evolves from (calcium driven) small oscillatory potentials SOP (on the left) to ``square-wave"-like bursting (on the right) across ``tapered"-like bursting (shown in the inset). Parameters values are provided in Section \ref{SEC: parameters}.}
\label{FIG 11}
\end{figure}

The single geometric attractor of (\ref{EQ: SIAM burst quantitative}) contains a continuum of different bursting wave forms modulated by the bifurcation and the unfolding parameters. Beyond the route to bursting studied in Section \ref{SEC: the routes},  Figure \ref{FIG 11} illustrates a situation where the bifurcation parameter and the affine unfolding parameters are fixed but where the two remaining unfolding parameters are modulated in a quasi static manner. The figure displays a variety of  waveforms that nevertheless share the same geometry of the bursting attractor as hysterethic paths in the universal unfolding of the winged cusp. For small autocatalytic feedback gain $k$, corresponding to low expression of fast sodium channels, the model emits small oscillatory potentials (SOP), on the left. Increasing this gain, the waveform smoothly evolves toward a classical `` square-wave" oscillation, on the right, after a transient ``tapered" bursting activity, shown in the inset (see \cite[page 376]{Izhikevich2007} and references therein for a discussion about the different bursting types). As in the case of the route from tonic spiking to bursting, the transition shown  in Fig. \ref{FIG 11} has  physiological relevance. For instance, a similar transition has been observed during development of neuronal cells \cite{Liu1998}.

\begin{figure}[h!]
\center
\includegraphics[width=0.9\textwidth]{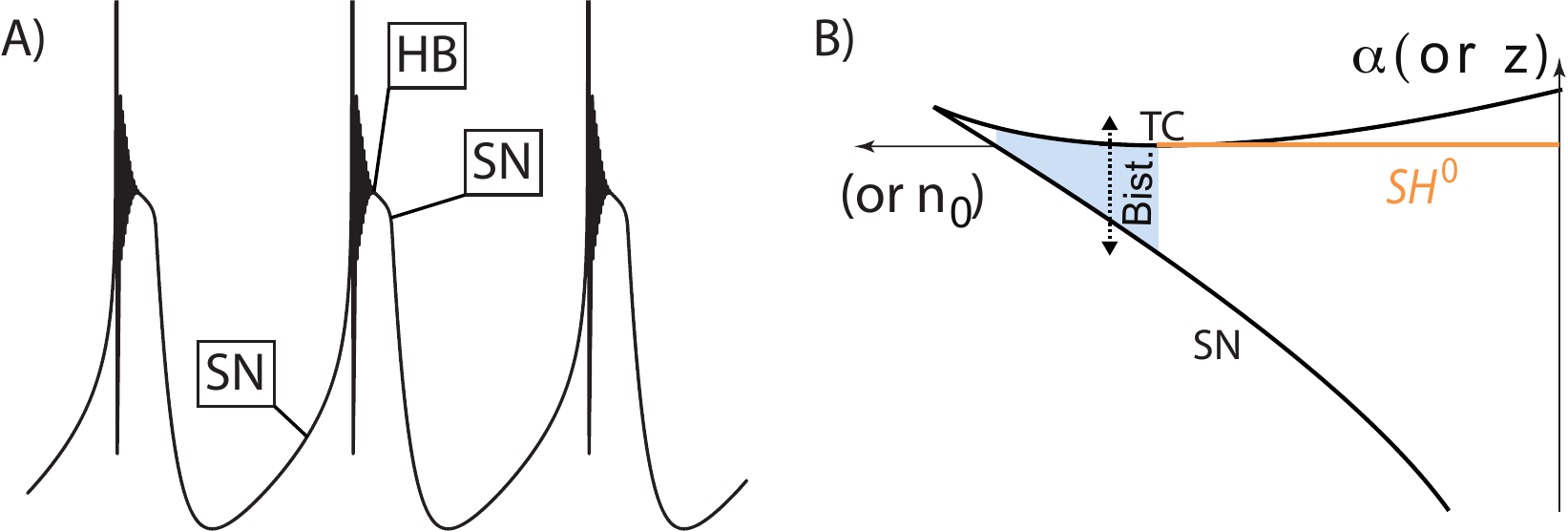}
\caption{{\bf Variant of the saddle-homoclinic bursting attractor in model (\ref{EQ: SIAM burst quantitative}).} {\bf A.} When the fast-slow subsystem (\ref{EQ: SIAM burst quantitative}a-\ref{EQ: SIAM burst quantitative}b) exhibits Type V excitability \cite{Franci2012}, the bistable range is of fold / fold type leading to a ``tapered" bursting waveform. Parameters values are provided in Section \ref{SEC: parameters}. {\bf B.} The hysteretic path associated to this type appears along path {\it ii)} of Fig. \ref{FIG 8}. {At the two ends of the bistable range, the up and down attractors are stable equilibria loosing stability in a saddle-node bifurcation}. Depending on the excitability subtype, the burst onset can either exhibit damped spiking oscillations ending in a Hopf bifurcation within the bistable range (a situation captured by the bifurcation diagram in \cite[Fig. 5.2]{Franci2012}) or a single action potential (a situation captured by the bifurcation diagram in \cite[Fig. 5.3]{Franci2012}).
}\label{FIG 12}
\end{figure}

The geometry of the ``tapered"-like bursting wave-form in Figure \ref{FIG 11} reveals  another subtlety of the winged cusp unfolding. In addition to broad regions of restorative and regenerative excitability, Fig.  \ref{FIG 8}A shows a small parametric region of mixed excitability (type V in the terminology of \cite{Franci2012}). Like regenerative phase portraits, phase portraits in this region have a persistent bistable range, but it is of fold/fold type, with a down-state that is a regenerative fixed point and a up-state that is either a restorative fixed point or a limit cycle (emerging from a Hopf bifurcation within or outside the bistable range). The bursting attractor observed in this region can be considered as a variant of the bursting attractor associated to regenerative excitability. Both bursting attractors share the same geometry of hysteretic paths in the unfolding of the winged cusp singularity but the fold/fold variant exhibits the peculiar wave form illustrated in Fig. \ref{FIG 12}, usually studied under the name of ``tapered" bursting in the literature, see e.g. \cite[page 376]{Izhikevich2007}.

It is remarkable that the four different wave forms shown in Figs. \ref{FIG 11}  and Fig. \ref{FIG 12} can be modeled by the same geometric attractor. A companion paper in preparation further investigates the physiological mechanisms that modulate the bursting wave within the unfolding of the winged cusp singularity.

\section{Conclusions}

The paper proposes that conductance based models exhibiting bursting attractors are organized by a winged cusp singularity. The geometry of the resulting attractor is classical (a hysteretic modulation of a slow-fast portrait over a rest-spike bistable range) but singularity theory is used to identify key parameters for the modulation of the bursting attractor.

The cusp singularity organizes the slow-fast phase portrait around the mirror hysteresis of Section \ref{SSEC: rest-spike in wcusp} in contrast to the standard hysteresis of classical phase portrait reductions of Hodgkin-Huxley model.

The bifurcation parameter has the convenient physiological interpretation of a ionic balance recently studied in \cite{Franci2013}. Its modulation through the transcritical variety of the cusp unfolding governs a geometric transition from tonic spiking to bursting in the three-timescale normal form (\ref{EQ: 3D cusp dynamics}): it provides a physiologically relevant route to bursting.

The affine unfolding parameter has the physiological interpretation of an ultraslow ionic current, typically driven by the intracellular calcium concentration. Its modulation provides the classical adaptation variable of the three time-scale bursting attractor.

The two remaining unfolding parameters have the physiological interpretation of a fast autocatalytic gain (the maximal sodium conductance) and of an average half activation potential, respectively. Their quasi static modulation evolves the bursting attractor across different bursting wave forms, consistently with what is observed experimentally in neuronal development, for instance.

In spite of the vast diversity of ion channels encountered in different neurons and the resulting vast diversity of regulation pathways, singularity theory and time scale separation suggest an apparent simplicity and universality in the underlying modulation mechanisms, as paths in the universal unfolding of the winged cusp. Those features are appealing to address system theoretic questions such as sensitivity, robustness, and homeostasis issues.

\section{Acknowledgments}

Prof. M. Golubitsky is gratefully acknowledged for insightful comments and suggestions during the visit of the first author at the Mathematical Bioscience Institute (Ohio State University).

\clearpage
\appendix

\section{Codimension 1 and 2 bifurcation varieties in (\ref{EQ: cusp dynamics})}
\label{SEC: trans variety}

The fixed point equation of (\ref{EQ: cusp dynamics}) is organized by a winged cusp at $x_{wcusp}:=\frac{1}{3}$, $\lambda_{wcusp}:=0$, $\alpha_{wcusp}:=-\frac{1}{27}$, $\beta_{wcusp}:=-\frac{1}{3}$, $\gamma_{wcusp}:=-2$. Codimension 1 transcritical and hysteresis bifurcation transition varieties in its unfolding are defined by
\begin{IEEEeqnarray}{rCl}\label{EQ: TC variety}
\alpha_{TC}(\beta,\gamma)&=&-\bar x_{TC}^3-(\bar\lambda_{TC}+\bar x_{TC})^2+\beta\bar x_{TC}-\gamma\bar x_{TC}(\bar\lambda_{TC}+\bar x_{TC})
\end{IEEEeqnarray}
with
\begin{IEEEeqnarray}{rCl}\label{EQ: TC variety 1}
\bar x_{TC}(\beta,\gamma)&=&\frac{\gamma^2-(\gamma^4+48\beta)^{1/2}}{12}\IEEEyessubnumber\\
\bar\lambda_{TC}(\beta,\gamma)&=&-\frac{\bar x_{TC}(2+\gamma)}{2}\IEEEyessubnumber
\end{IEEEeqnarray}
and
\begin{IEEEeqnarray}{rCl}\label{EQ: HY variety}
\alpha_{HY}(\beta,\gamma)&=&-\bar x_{HY}^3-(\bar\lambda_{HY}+\bar x_{HY})^2+\beta\bar x_{HY}-\gamma\bar x_{HY}(\bar\lambda_{HY}+\bar x_{HY})
\end{IEEEeqnarray}
with
\begin{IEEEeqnarray}{rCl}\label{EQ: HY variety 1}
\bar x_{HY}(\gamma)&=&-\frac{1+\gamma}{3}\IEEEyessubnumber\\
\bar\lambda_{HY}(\beta,\gamma)&=&\frac{\beta-3\bar x^2_{HY}-\bar x_{HY}(2+2\gamma)}{2+\gamma}\IEEEyessubnumber
\end{IEEEeqnarray}
respectively.

The codimension 2 pitchfork variety is defined by
\begin{IEEEeqnarray}{rCl}\label{EQ: PF variety}
\gamma_{PF}(\beta)&=&\left(3\beta + \left(9\beta^2 + \frac{1}{27}\right)^{1/2}\right)^{1/3} - \frac{1}{\left(3\beta + \left(9\beta^2 + \frac{1}{27}\right)^{1/2}\right)^{1/3}} - 1\IEEEyessubnumber\\
\alpha_{PF}(\beta)&=&-\bar x_{PF}^3-(\bar\lambda_{PF}+\bar x_{PF})^2+\beta\bar x_{PF}-\gamma\bar x_{PF}(\bar\lambda_{PF}+\bar x_{PF})
\end{IEEEeqnarray}
with
\begin{IEEEeqnarray}{rCl}\label{EQ: PF variety 1}
\bar x_{PF}(\beta)&=&-\frac{1+\gamma_{PF}}{3}\IEEEyessubnumber\\
\bar\lambda_{PF}(\beta)&=&\frac{\beta-3\bar x_{PF}^2-\bar x_{PF}(2+2\gamma_{PF})}{2+\gamma_{PF}}\IEEEyessubnumber
\end{IEEEeqnarray}

\section{Proofs}

\subsection{Proof of Theorem \ref{THM: cusp bist}}
\label{SSEC: cusp bist proof}

We rely on geometric singular perturbation arguments \cite{Fenichel1979,Jones1995,Krupa2001c,Krupa2001b,Krupa2001a}. The reduced dynamics associated to (\ref{EQ: cusp dynamics}), evolving on the slow time scale $\tau=\varepsilon t$, is given by
\begin{IEEEeqnarray}{rCl}\label{EQ: cusp dynamics reduced}
0&=&G_{wcusp}^s(x,\lambda+y;\ \alpha,\beta,\gamma)\IEEEyessubnumber\\
\dot y&=&x-y,\IEEEyessubnumber
\end{IEEEeqnarray}
whereas the associated layer dynamics, evolving on the fast time scale $t$, is given by
\begin{IEEEeqnarray}{rCl}\label{EQ: cusp dynamics layer}
\dot x&=&G_{wcusp}^s(x,\lambda+y;\ \alpha,\beta,\gamma)\IEEEyessubnumber\\
\dot y&=&0.\IEEEyessubnumber
\end{IEEEeqnarray}

\begin{figure}[h!]
\center
\includegraphics[width=0.85\textwidth]{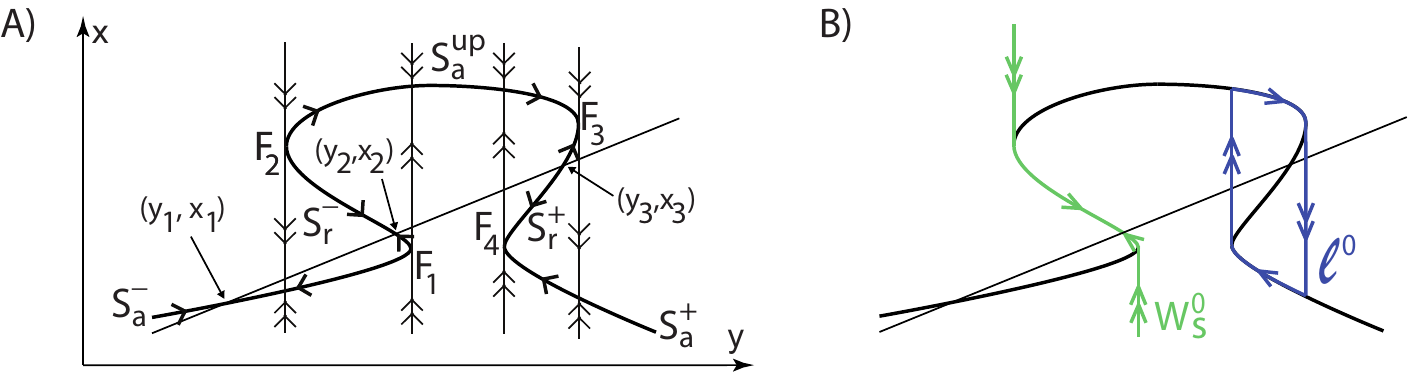}
\caption{{\bf Slow fast dynamics associated to (\ref{EQ: cusp dynamics reduced}-\ref{EQ: cusp dynamics layer}).}}\label{SUPP FIG 1}
\end{figure}

We construct the singular bistable phase portrait starting from the degenerate situation in Fig. \ref{FIG 3} center, corresponding to a pitchfork bifurcation. The same qualitative phase portrait is obtained on the pitchfork variety (\ref{EQ: PF variety}) for all $\beta>\beta_{wcusp}$. Perturbing $\gamma$ out of the pitchfork variety, but remaining on the transcritical variety defined by (\ref{EQ: TC variety}), the phase portrait perturbs to one of the two qualitative situations in Fig. \ref{FIG 3} center - top or bottom. Finally, for $\lambda$ below and sufficiently near $\lambda_{TC}(\beta,\gamma)$ and $\alpha$ below and sufficiently near $\alpha_{TC}(\beta,\gamma)$ one obtains the qualitative slow-fast dynamics in Fig. \ref{SUPP FIG 1}A, which leads to the singular phase-portrait in Fig. \ref{SUPP FIG 1}B. The following lemma summarizes this construction.

\begin{lemma}\label{LEM: cusp bist 1}
For all $\beta>\beta_{wcusp}$, there exists $\Delta_\gamma>0$ such that, for all $\gamma\in(\gamma_{PF}(\beta)-\Delta_\gamma,\gamma_{PF}(\beta)+\Delta_\gamma)$, there exists $\Delta_\lambda>0$ such that, for all $\lambda\in(\bar\lambda_{TC}(\beta,\gamma)-\Delta_\lambda,\bar\lambda_{TC}(\beta,\gamma))$, there exists $\Delta_\alpha>0$ such that, for all $\alpha\in(\alpha_{TC}(\beta,\gamma)-\Delta_\alpha,\alpha_{TC}(\beta,\gamma))$, the following hold  {\rm (refer to Fig. \ref{SUPP FIG 1} - left for the notation)}:
\begin{itemize}
\item[{\it i)}] The critical manifold of the slow-fast dynamics (\ref{EQ: cusp dynamics reduced}-\ref{EQ: cusp dynamics layer}) has a mirrored hysteresis shape. In particular, it is composed of the attractive branches $S_a^-$, $S_a^+$, and $S_a^{up}$, the repelling branches $S_r^-$ and $S_r^+$, and the four folds $F_i$, $i=1,\ldots,4$, connecting them.
\item[{\it ii)}] There are exactly three nullcline intersection $(y_i,x_i)$, $i=1,\ldots,3$, belonging to $S_a^-$, $S_r^-$, and $S_r^+$, respectively.
\end{itemize}
\end{lemma}

A direct geometric inspection reveals the presence of a singular periodic orbit $\ell^0$ and a singular saddle separatrix $W_s^0$. These objects persist for $\epsilon>0$, as proved in the following lemma, which proves Theorem \ref{THM: cusp bist}.

\begin{lemma}\label{LEM: cusp bist 2}
Let $(y_i,x_i)$, $i=1,\ldots,3$ be defined as in the statement of Lemma \ref{LEM: cusp bist 1} {\it ii)}. For all $\lambda,\alpha,\beta,\gamma$ satisfying conditions of Lemma \ref{LEM: cusp bist 1}, there exists $\bar\varepsilon$ such that, for all $\varepsilon\in (0,\bar\varepsilon)$

{\it i)} $(y_1,x_1)$ is locally exponentially stable, $(y_2,x_2)$ is a hyperbolic saddle, and $(y_3,x_3)$ is locally exponentially unstable.

{\it ii)} There exists an exponentially stable relaxation oscillation limit cycle $\ell^\varepsilon$ surrounding $(y_3,x_3)$.

{\it iii)} The stable manifold $W_s^\varepsilon$ of $(y_2,x_2)$ separates the basin of attraction of $(y_1,x_1)$ and $\ell^\varepsilon$.

\end{lemma}
{\bf Proof of Lemma \ref{LEM: cusp bist 2}.}

{\it i)}  From Lemma \ref{LEM: cusp bist 1}, the fixed point $(y_1,x_1)$ belongs to the attractive branch $\mathcal S_a^-$ of the critical manifold $\mathcal S$. Moreover, it is an exponentially stable fixed point of the the reduced dynamics (\ref{EQ: cusp dynamics reduced}). From standard persistence arguments \cite{Fenichel1979}, there exists $\bar\varepsilon_1$ such that, for all $\epsilon\in(0,\bar\varepsilon_1]$, $(y_1,x_1)$ is an exponentially stable fixed point of (\ref{EQ: cusp dynamics}). The fixed point $(y_2,x_2)$ belongs to the repelling branch $\mathcal S_r^-$ of the critical manifold $\mathcal S$. Moreover, it is an exponentially stable fixed point of the reduced dynamics (\ref{EQ: cusp dynamics reduced}). Again from \cite{Fenichel1979}, there exists $\bar\varepsilon_2$ such that, for all $\epsilon\in(0,\bar\varepsilon_2]$, there exists an exponentially unstable local invariant manifold $W_{s,loc}^{\varepsilon}$ such that all trajectories starting in $W_{s,loc}^{\varepsilon}$ approach $(y_2,x_2)$ exponentially fast. $W_{s,loc}^{\varepsilon}$ is the local stable manifold of $(y_2,x_2)$. Its unstable manifold is given by the fiber of the unstable manifold of $W_{s,loc}^{\varepsilon}$ passing through $(y_2,x_2)$. The fixed point $(y_3,x_3)$ belong to the repelling branch $\mathcal S_r^+$ of the critical manifold $\mathcal S$, moreover it is an exponentially unstable fixed point of the reduced dynamics (\ref{EQ: cusp dynamics reduced}). By \cite{Fenichel1979}, there exists $\bar\varepsilon_3>0$ such that, for all $\epsilon\in(0,\bar\varepsilon_3]$, $(y_3,x_3)$ is an exponentially unstable fixed point of (\ref{EQ: cusp dynamics}).

{\it ii)} The slow fast dynamics possesses a singular periodic orbit $\ell^0$ (See Fig. \ref{SUPP FIG 1}). Following \cite{Krupa2001c}, there exists $\bar\varepsilon_4$ such that, for all $\epsilon\in(0,\bar\varepsilon_4]$, there exists an exponentially stable relaxation oscillation limit cycle $\ell^\varepsilon$ surrounding $(y_3,x_3)$.

{\it iii)} In backward time, trajectories of the reduced dynamics (\ref{EQ: cusp dynamics reduced}) starting on $\mathcal S_r^-$ in a neighborhood of $(y_2,x_2)$ approach either the fold $\mathcal F_1$ or the fold $\mathcal F_2$. Following \cite{Krupa2001b}, there exists $\bar\varepsilon_5$ such that, for all $\epsilon\in(0,\bar\varepsilon_5]$, all trajectories starting in the local stable manifold $W_{s,loc}^{\varepsilon}$ approach (in backward time) either the fold $\mathcal F_1$ or the fold $\mathcal F_2$ along an invariant manifold $W_s^\epsilon$, which continues after the fold singularities roughly parallel to trajectories of the layer problem. Therefore, the branch that continues after $\mathcal F_1$ extends to $x=-\infty$, whereas the branch that continues after $\mathcal F_2$ extends to $x=+\infty$. The invariant manifold $W_s^\epsilon$ is the saddle stable manifold and separates the plane in two disconnected regions that contain, respectively, the two attractors $(y_1,x_1)$ and $\ell^\varepsilon$.

Items {\it i)}, {\it ii)}, and {\it iii)} are proved by picking $\bar\varepsilon=\min_{i=1,\ldots,5}\bar\varepsilon_i$.
\hfill$\square$\\

\hfill$\square$\\

\subsection{Proof of Theorem \ref{THM: cusp burst}}
\label{SSUB: burst att proof}

\begin{figure}[h!]
\center
\includegraphics[width=0.85\textwidth]{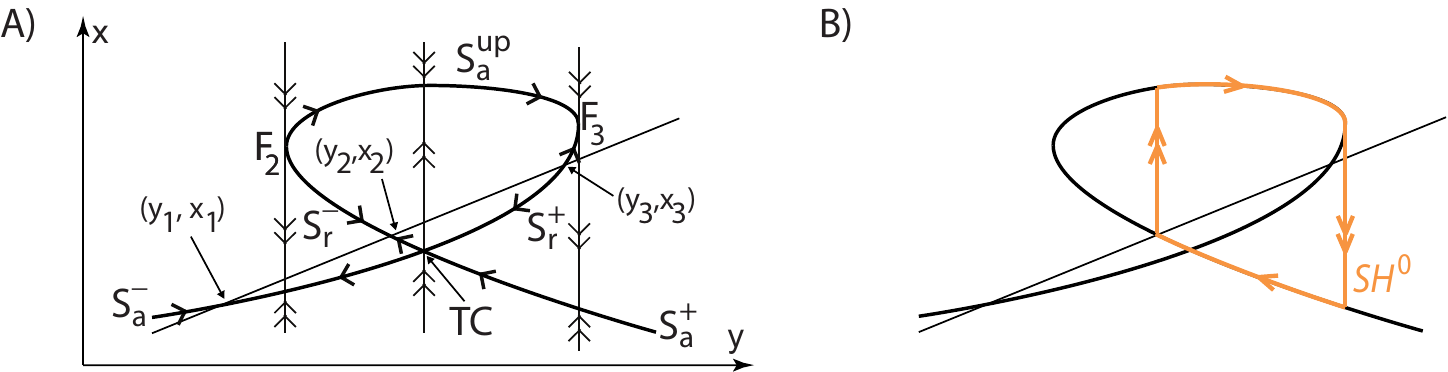}
\caption{{\bf Slow fast dynamics associated to (\ref{EQ: cusp dynamics reduced}-\ref{EQ: cusp dynamics layer}).}}\label{SUPP FIG 2}
\end{figure}

\begin{figure}[h!]
\center
\includegraphics[width=0.85\textwidth]{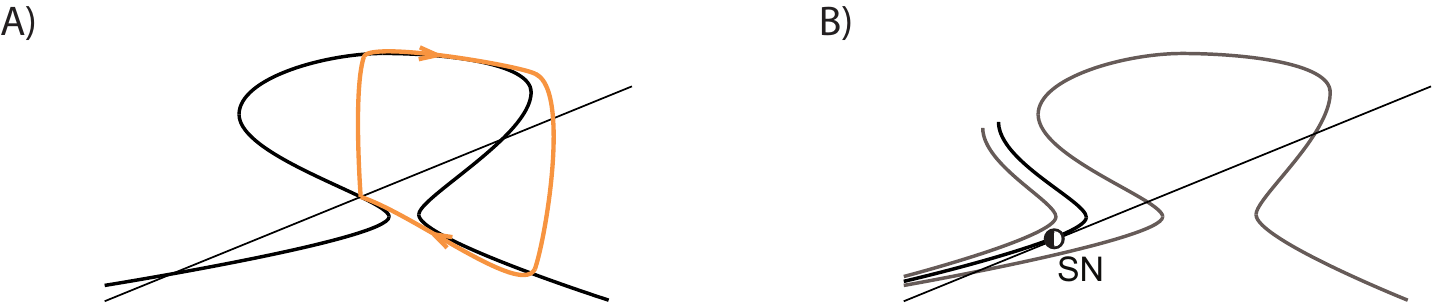}
\caption{{\bf Phase portrait of (\ref{EQ: cusp dynamics}), with parameters and the function $y_\infty$ satisfying conditions of Lemma \ref{LEM: cusp SH}}. {\bf A.} Singularly perturbed saddle-homoclinic trajectory for $\alpha=\alpha_{SH}^0+\alpha_c(\sqrt{\varepsilon})$. {\bf B.} Saddle-node bifurcation.
%The three plots of the function $f_2$, defined in the proof of Lemma \ref{LEM: cusp SH}, correspond to different values of $\alpha$: 1, $\alpha_{SN}<\alpha<\alpha_{SH}^0$; 2, $\alpha=\alpha_{SN}$, 3: $\alpha<\alpha_{SN}$.
}\label{SUPP FIG 3}
\end{figure}

Starting from a set of parameter satisfying the condition of Lemma \ref{LEM: cusp bist 1} and increasing $\alpha$ to $\alpha=\alpha_{TC}(\beta,\gamma)$ the two folds $F_1$ and $F_4$ in Fig. \ref{SUPP FIG 1}A approach each other and eventually collide in a transcritical singularity $TC$, as in the slow-fast dynamics in Figure \ref{SUPP FIG 2}A. A direct geometrical inspection reveals the presence of a singular saddle-homoclinic trajectory $SH^0$ (Fig. \ref{SUPP FIG 2}B) for which the transcritical singularity serves as connecting point. This homoclinic orbit persists for $\varepsilon>0$, as sketched in Figure \ref{SUPP FIG 3}A. On the contrary, decreasing $\alpha$ the two folds move away from each other until the left branch of the mirrored hysteresis is tangent to the $y$ nullcline at a saddle-node bifurcation $SN$ and eventually remains on its left, as in Fig. \ref{SUPP FIG 3}B. The following lemma summarizes this analysis. For its statement, we refer to Figures \ref{SUPP FIG 1} and \ref{SUPP FIG 2}.

\begin{lemma}\label{LEM: cusp SH}
For all $\beta>\beta_{wcusp}$, there exists $\Delta_\gamma>0$ such that, for all $\gamma\in(\gamma_{PF}(\beta)-\Delta_\gamma,\gamma_{PF}(\beta)+\Delta_\gamma)$, there exist $\Delta_\lambda>0$ such that, for all $\lambda\in(\bar\lambda_{TC}(\beta,\gamma)-\Delta_\lambda,\bar\lambda_{TC}(\beta,\gamma))$, {there exists $\bar\varepsilon>0$, such that, for $\varepsilon\:in(0,\bar\varepsilon]$ the following hold}:

{\it i)} Let $\alpha_{SH}^0:=\alpha_{TC}(\beta,\gamma)$.  There exists a smooth function $\alpha_c(\cdot)$ defined on $\left[0,\sqrt{\bar\varepsilon}\right]$, and satisfying $\alpha_c(0)=0$ and $\alpha_c(\sqrt{\varepsilon})<0$, for all $\varepsilon\in(0,\bar\varepsilon]$, such that, for $\alpha=\alpha_{SH}^\varepsilon:=\alpha_{SH}^0+\alpha_c(\sqrt{\varepsilon})$, (\ref{EQ: cusp dynamics}) has an unstable saddle-homoclinic orbit $SH^\varepsilon$.

{\it ii)} For all $\varepsilon>0$, there exists $\alpha_{SN}<\alpha_{SH}^0$ such that (\ref{EQ: cusp dynamics}) has a non-degenerate saddle-node bifurcation for $\alpha=\alpha_{SN}$ at which the node $(y_1,x_1)$ and the saddle $(y_2,x_2)$ merge.

{\it iii)} For all $\alpha\in[\alpha_{SN},\alpha_{SH}^0)$ the nullcline intersection $(y_3,x_3)$ belongs to the repelling branch $S_r^+$ (where $S_r^+$ is defined as in Figure \ref{SUPP FIG 1}A).

{\it iv)} {For all $\alpha\in[\alpha_{SN},\alpha_{SH}^\varepsilon]$, there exists an exponentially stable relaxation oscillation limit cycle $\ell^\varepsilon$ surrounding $(y_3,x_3)$.}

{\it v)} {There exists $\alpha_{FLC}^\varepsilon\in(\alpha_{SH}^\varepsilon,\alpha_{SH}^0)$ such that, for $\alpha=\alpha_{FLC}^\varepsilon$, the family $P^\varepsilon$ of stable periodic orbits merge at a fold limit cycle bifurcation with the family of unstable periodic orbits $Q^\varepsilon$ emerging from the unstable saddle-homoclinic bifurcation.}
\end{lemma}

{\bf Proof of Lemma \ref{LEM: cusp SH}.}
% Again, we start by considering the degenerate situation $\alpha=\alpha_{PF}(\beta),\gamma=\gamma_{PF}(\beta)$ and $\lambda=\bar\lambda_{PF}(\beta)$, for which one obtains the phase portrait in Figure \ref{FIG 3} center.
{\it i)} For $\gamma$ in a neighborhood of  $\gamma_{PF}(\beta)$, $\alpha=\alpha_{TC}(\beta,\gamma)$ and $\lambda$ smaller than and sufficiently near to $\bar\lambda_{TC}(\beta,\gamma)$, there are exactly three nullcline intersections $(y_i,x_i)$, $i=1,2,3$, belonging to the attractive branch $S_a^-$, the repelling branch $S_r^-$, and the repelling branch $S_r^+$, respectively.  {Relying on the results in} \cite{Krupa2001a} {and following exactly the same steps as} \cite[Section 6.1]{Franci2012}{, we can find $\varepsilon_1$, such that the existence part of the point {\it i)} holds with $\bar\varepsilon=\varepsilon_1$}. The resulting saddle-homoclinic trajectory is sketched in Fig. \ref{SUPP FIG 3}A. To prove that such homoclinic trajectory is unstable, recall that the stability of a saddle-homoclinic orbit is determined by the saddle quantity $\sigma$, that is, the trace of the Jacobian computed at the saddle and at the saddle-homoclinic bifurcation: if $\sigma>0$ (resp. $\sigma<0$) the homoclinic orbit is unstable (resp. stable). The Jacobian $J_{SH}$ of (\ref{EQ: cusp dynamics}) computed at $(y_2,x_2)$ at the saddle-homoclinic bifurcation has the form
\begin{equation*}
J_{SH}=\left(\begin{array}{cc}
a & b\\ \varepsilon c & -\varepsilon d
\end{array}
\right),\quad a,d>0,\ b,c\in\mathbb R.
\end{equation*}
Therefore the saddle quantity $\sigma=a-\varepsilon d>0$, for all $0<\varepsilon<a/d$.\\

 {\it ii)} For $\gamma$ in a neighborhood of  $\gamma_{PF}(\beta)$, $\alpha=\alpha_{TC}(\beta,\gamma)$, and $\lambda$ smaller than and sufficiently near to $\bar\lambda_{TC}(\beta,\gamma)$, the (cubic) fixed point equation $G_{wcusp}(x,\lambda+x;\alpha,\beta,\gamma)$ has three roots, corresponding to the three fixed point $(y_i,x_i)$, $i=1,2,3$ of point {\it i)}. Decreasing $\alpha$, the two smaller roots (corresponding to the fixed point $(y_1,x_1)$ and $(y_2,x_2)$) approach each other and eventually merge in a quadratic zero for $\alpha=\alpha_{SN}$ corresponding to a non-degenerate saddle-node bifurcation.

{\it iii)} We prove the statement for $\gamma=\gamma_{PF}(\beta)$ since, by continuity, the same will hold in a neighborhood. When $\gamma=\gamma_{PF}(\beta)$, $\alpha=\alpha_{PF}(\beta)$, and $\lambda$ is smaller than and sufficiently near to $\bar\lambda_{PF}(\beta)$, the nullcline intersection $(y_3,x_3)$ lies on $S^+_r$. By continuity, the same is true for all $\alpha$ close to $\alpha_{PF}(\beta)$. Since the value $\alpha_{SN}\uparrow\alpha_{PF}(\beta)$ continuously as $\lambda\uparrow\bar\lambda_{PF}(\beta)$, one can pick $\lambda$ sufficiently close to $\bar\lambda_{PF}(\beta)$ such that $(y_3,x_3)$ lies on $S^+_r$ for all $\alpha\in[\alpha_{SN},\alpha_{PF}(\beta))$.

{\it iv)} {By points {\it ii)} and {\it iii)} above and the same arguments as the proof of point {\it ii)} in Lemma} \ref{LEM: cusp SH}{, we can find $\varepsilon_2>0$ such that, for all $\varepsilon\in(0,\min(\varepsilon_1,\varepsilon_2)$, where $\varepsilon_1$ is defined as in the proof of point {\it i)} above, and all $\alpha\in[\alpha_{SN},\alpha_{SH}^\varepsilon)$, there exists a periodic orbit $\ell^\varepsilon$ surrounding $(y_3,x_3)$ and, moreover, this periodic orbit is exponentially stable. For $\alpha=\alpha_{SH}^\varepsilon$, the stable periodic orbit co-exist with the unstable homoclinic orbit, since by} \cite[Theorem 3.5]{Chow1994}{, a branch of stable periodic orbits cannot end at an unstable homoclinic bifurcation.}

{\it v)} {The existence of $\alpha_{FLC}^\varepsilon\in(\alpha_{SH}^\varepsilon,\alpha_{SH}^0)$ satisfying the statement follows by two main observations. First, again by} \cite[Theorem 3.5]{Chow1994} {there exists a family $Q^\varepsilon$ of unstable periodic orbits emergenging at $\alpha=\alpha^\varepsilon_{SH}$ from the unstable homoclinic bifurcation. Second, simple geometric arguments show that for $\alpha=\alpha_{SH}^0$ (and $\varepsilon$ sufficiently small) no periodic orbit can exists. The existence of the fold limit cycle bifurcation then follows by noticing that the fold limit cycle is the only planar bifurcation of periodic orbits not involving a Hopf point and that both the unstable homoclinic bifurcation and the fold limit cycle bifurcation are generically found in the unfolding of the degenerate situation in which the saddle quantity $\sigma$ is zero, corresponding to a resonant homoclinic orbit. The unfolding of this bifurcation, also called resonant side-switching, is detailed in} \cite[Theorem A]{Chow1990}.
\hfill$\square$\\

\begin{figure}[h!]
\center
\includegraphics[width=0.75\textwidth]{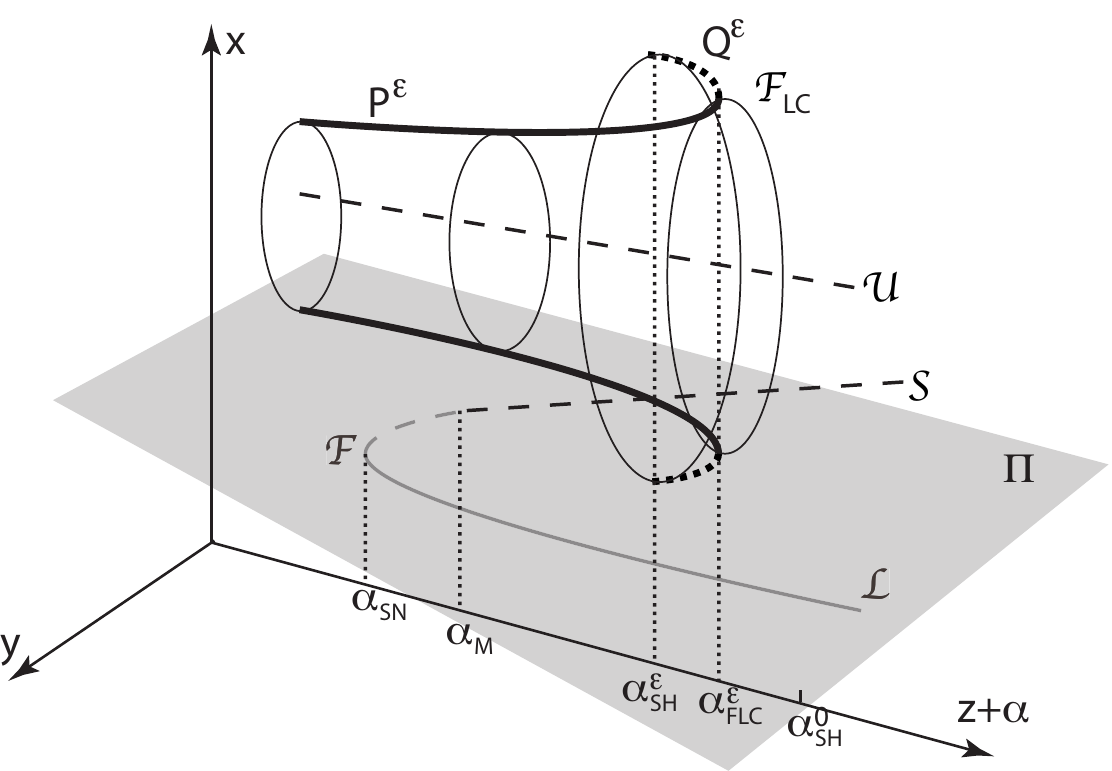}
\caption{{\bf Bifurcation diagram of (\ref{EQ: 3D cusp dynamics}a-\ref{EQ: 3D cusp dynamics}b) with respect to $z$ and parameters satisfying condition of Lemma \ref{LEM: cusp SH}}. See the main text describing Fig. \ref{FIG 5} for the notation.}\label{SUPP FIG 4}
\end{figure}

Figure \ref{SUPP FIG 4} summarizes the results in Lemma \ref{LEM: cusp SH}.

%Lemmas \ref{LEM: cusp bist 2} and \ref{LEM: cusp SH} imply that there exists $\bar\varepsilon$ such that, for all $\varepsilon\in(0,\bar\varepsilon]$ and parameters satisfying condition of Lemma \ref{LEM: cusp SH}, the qualitative bifurcation diagram of (\ref{EQ: 3D cusp dynamics}a-\ref{EQ: 3D cusp dynamics}b) with respect to $z$ is as in Figure \ref{SUPP FIG 4}. {Let us develop. Nearby the unstable homoclinic bifurcation there is a saddle-node limit cycle bifurcation $\alpha_{FLC}$ in which the stable limit cycle merges with the unstable limit cycle born at the unstable homoclinic bifurcation. Both bifurcations arise in the unfolding of the degenerate situation in which the saddle quantity $\sigma$ is zero, corresponding to a resonant homoclinic orbit. The unfolding of this bifurcation, also called resonant side-switching, is detailed in} \cite[Theorem A]{Chow1990}. The qualitative bifurcation diagram in Fig. \ref{SUPP FIG 4} follows.

We now follow \cite{Terman1991,Su2004} to derive suitable conditions on the four parameters $a,b,c$ in (\ref{EQ: 3D cusp dynamics}c) such that $z$ hysteretically modulates (\ref{EQ: 3D cusp dynamics}a-\ref{EQ: 3D cusp dynamics}b) along its rest spike bistable range. To this aim, note that the minimum value of $y$ along the family of singular periodic orbits $\ell^0$ and the singular homoclinic trajectory $SH^0$ for $\alpha\in[\alpha_{SN},\alpha{SH}^0]$ (see Figure \ref{FIG 4}) is necessarily strictly larger then the maximum value of $y$ along the branch of stable fixed points $\mathcal L$ and at the saddle-node bifurcation $\mathcal F$. By persistence arguments, the same holds true in the nonsingular case for $\varepsilon$ sufficiently small. It follows that there exists a plane $\Pi$ in the three-dimensional space $x,y,z$ that, for $\alpha\in[\alpha_{SN},\alpha_{SH}^0]$, never intersects the family of stable periodic orbits $P^\varepsilon$ and the branch of stable fixed points $\mathcal L$, and that intersects once the branch of saddle points $\mathcal S$, say, for $\alpha=\alpha_M$  (see Fig. \ref{SUPP FIG 4}). Clearly, $\Pi$ splits $\mathbb R^3$ in two open half-spaces. Let $\Pi^c_\ell$ be the half space containing the family of singular periodic orbits. Then we pick $a,b,c$ such that $\Pi=\{(x,y,z):\ -z+ax+by+c=0\}$ and $-z+ax+by+c>0$, for $(x,y,z)\in\Pi^c_\ell$. Under these conditions on $a,b,c$, Theorem \ref{THM: cusp burst} follows along the same line as the proofs in \cite{Terman1991} (for the analysis near the branch of the stable steady states and the ``jump up" at the fold bifurcation) and \cite{Su2004} (for the analysis near the branch of periodic orbits and the ``jump down" at the fold limit cycle bifurcation). \hfill$\square$

\section{Parameter for numerical simulations in Figs. \ref{FIG 6}, \ref{FIG 7}, \ref{FIG 11}, \ref{FIG 12}}
\label{SEC: parameters}

For the sake of an easy numerical implementation and the reproduction of ``nice" {time series}, we suggest the following piecewise linear approximation of (\ref{EQ: SIAM burst quantitative})
\begin{IEEEeqnarray*}{rCl}
\dot V&=&kV-\frac{V^3}{3}-(n+n_0)^2+I-z\IEEEyessubnumber\\
\dot n&=&\varepsilon_n\left(\hat n_\infty(V-V_0) -n \right) \IEEEyessubnumber\\
\dot z&=&\varepsilon_z(\hat z_\infty(V-V_1)-z)\IEEEyessubnumber
\end{IEEEeqnarray*}
where
\begin{IEEEeqnarray*}{rCl}
\hat n_\infty(V-V_0)&:=&\begin{cases}
k^n_-(V-V_0) & \text{if } V < V_0,\\
k^n_+(V-V_0) & \text{if } V \geq V_0.
\end{cases}
\end{IEEEeqnarray*}
with $0\leq k^n_-< 1$ and $k^n_+> 1$, and
\begin{IEEEeqnarray*}{rCl}
\hat z_\infty(V-V_1)&:=&\begin{cases}
k^z_-(V-V_1) & \text{if } V < V_1,\\
k^z_+(V-V_1) & \text{if } V \geq V_1.
\end{cases}
\end{IEEEeqnarray*}
with $0\leq k^z_-<1$ and $k^z_+>1$.

Parameters used in Figs. \ref{FIG 6} and \ref{FIG 7} are $k=1$, $I=11/3$, $\varepsilon_n=0.02$, $\varepsilon_z=0.0005$, $V_0=-0.5$, $k^n_-=0.4$, $k^n_+=7$, $V_1=-1$, $k^z_-=0$, $k^z_+=50$. The bifurcation parameter is $n_0=0.3$ in Figure \ref{FIG 6} left and $n_0=-1.1$ in Fig. \ref{FIG 6} right. In Fig. \ref{FIG 7} $n_0$ is linearly (in time) decreased from $0.3$ to $-1.1$.

Parameters used in Figs. \ref{FIG 11}A are $n_0=-1.1$, $I=11/3$, $\varepsilon_n=0.02$, $\varepsilon_z=0.0001$, $k^n_-=0.9$, $k^n_+=7$, $V_1=-1.2$, $k^z_-=0$, $k^z_+=100$. The time-varying parameters $k$ and $V_0$ evolve as $k(t)=-0.5+2.5t/T$ and $V_0(t)=-0.5-0.75\ \min(1,\ 1.3t/T)$.

%Parameters used in Figs. \ref{FIG 11}A are $n_0=-1.1$, $k=1$, $I=11/3$, $\varepsilon_n=0.02$, $\varepsilon_z=0.0005$, $V_0=-0.5$, $k^n_-=0.4$, $k^n_+=7$, $V_1=-1$, $k^z_-=0$, $k^z_+=50$.
%
%Parameters used in Figs. \ref{FIG 11}B are $n_0=-0.2$, $k=1$, $I=11/3$, $\varepsilon_n=0.02$, $\varepsilon_z=0.0002$, $V_0=-0.5$, $k^n_-=0.05$, $k^n_+=7$, $V_1=-0.5$, $k^z_-=0$, $k^z_+=50$.
%
%Parameters used in Figs. \ref{FIG 11}C are $n_0=-0.2$, $k=0.12$, $I=11/3$, $\varepsilon_n=0.02$, $\varepsilon_z=0.001$, $V_0=-0.5$, $k^n_-=0.4$, $k^n_+=7$, $V_1=-0.5$, $k^z_-=0$, $k^z_+=50$.

Parameters used in Figs. \ref{FIG 12}A are $n_0=-0.2$, $k=0.7$, $I=11/3$, $\varepsilon_n=0.02$, $\varepsilon_z=0.001$, $V_0=-1.25$, $k^n_-=0.9$, $k^n_+=7$, $V_1=-1.2$, $k^z_-=0$, $k^z_+=50$.

\section{Parameter for numerical simulations of the Hodgkin-Huxley model in Section \ref{SSSE: winged-cusp phase plane in HHCa}}
\label{SEC: parameters HH}
All the parameter and activation and inactivation rates are taken from the original paper \cite{Hodgkin1952}. The time constants $\tau_x(V)$ and steady state functions $x_\infty(V)$ are related to the activation and inactivation rates $\alpha_x(V)$ and $\beta_x(V)$, $x=m,n,h$, as follows
$$\tau_x(V)=\frac{1}{\alpha_x(V)+\beta_x(V)},\quad x_\infty(V)=\frac{\alpha_x(V)}{\alpha_x(V)+\beta_x(V)}.$$
%In accordance with the original paper \cite{Hodgkin1952}, the potassium half-activation $V_{1/2,K}=10mV$.

\bibliographystyle{unsrt}
\bibliography{../../../../../../Dropbox/Orchestron_bibliography/orchestron}

\end{document}